
\input amssym.def \input amssym 

\magnification=1200
\hsize=160truemm
\vsize=240truemm
\parskip=\smallskipamount

\def\qed{\hfill$\qquad\square$\par\medbreak}
\def\bmatrix#1{\left| \matrix{#1} \right|}

\def\wt #1{\widetilde{#1}}
\def\wl #1{\overline{#1}}
\def\sier#1{{\cal O}_{#1}}
\def\ttw(#1){T^2_X(#1)}
\def\tf{{\tilde f}}

\def\C{\Bbb C}
\def\P{{\Bbb P}}

\def\so{{\cal O}}
\def\cc{{\cal C}}
\def\ciz{{\cal I}_Z}
\def\al{\alpha}
\def\be{\beta}
\def\ga{\gamma}
\def\la{\lambda}
\def\de{\delta}
\def\ze{\zeta}
\def\om{\omega}
\def\vp{\ph}
\def\frac#1#2{{{#1}\over{#2}}}
\def\ip<#1,#2>{\langle#1,#2\rangle}
\def\Cork{\mathop{\rm Cork}}
\def\Hom{\mathop{\rm Hom}}
\def\Rank{\mathop{\rm Rank}}
\def\lra{\longrightarrow}
\def\ep{\varepsilon}
\def\ph{\varphi}

\def\cn{\mathbin{:}}
\long\def\comment#1\endcomment{}
\def\cite#1{[#1]}
\def\roep #1.{\medbreak\noindent{\sl #1\/}.\enspace}
\def\endroep{\par\medbreak}

\def\kl(#1){{\scriptscriptstyle(#1)}}

\newcount\secno
\newcount\thmno
\outer\def\section#1\par{\advance\secno by 1 \thmno=0 
  \vskip0pt plus.3\vsize\penalty-250
  \vskip0pt plus-.3\vsize\bigskip\vskip\parskip
  \message{#1}\leftline{\bf\the\secno.\enspace#1}
  \nobreak\smallskip\noindent}
\def\num{\global\advance\thmno by 1
\the\secno.\the\thmno}

\def\rfbs{Bayer--Stillman}
\def\rfds{Drewes--Stevens 1996}
\def\rfbra{Brawner 1997}
\def\rfbrap{Brawner 1996} 
\def\rfcm{Ciliberto--Miranda  1992}
\def\rfcp{Ciliberto--Pareschi  1995}
\def\rfduv{Du Val 1933}
\def\rfdonmor{Donagi--Morrison 1989}
\def\rfep{Epema 1983}
\def\rfrch{Reid 1997}
\def\rfr{Reid 1989}
\def\rfwasq{Wahl 1997}
\def\rfwacy{Wahl 1998}
\def\rfrohn{Rohn 1884}
\def\rfschr{Schreyer 1986}
\def\rfsth{Stevens 1996}  

{\font\largebf=cmbx10 scaled \magstep2
\centerline{\largebf Rolling Factors Deformations}
\medskip
\centerline{\largebf and Extensions of Canonical Curves}
\bigskip

\centerline{\bf Jan Stevens}}
\bigskip

\bigskip

\bigskip

An easy dimension count shows that not all canonical curves are 
hyperplane sections of $K3$ surfaces. A surface with a given curve
as hyperplane section is called an extension of the curve.
With this terminology, the general canonical curve has only trivial 
extensions, obtained by taking a cone over the curve.
In this paper we concentrate on extensions of tetragonal curves.

The extension problem  is related to deformation theory for cones.
This is best seen in terms of equations.
Suppose we have coordinates $(x_0\cn\cdots\cn x_n\cn t)$ on
$\P^{n+1}$ with the special hyperplane section  given by $t=0$.
We describe an extension $W$ of a variety $V\colon f_j(x_i)=0$ 
by a system of equations $F_j(x_i,t)=0$ with $F_j(x_i,0)=f_j(x_i)$. 
We write $F_j(x_i,t)=f_j(x_i)+tf_j'(x_i)
+\cdots+a_jt^{d_j}$, where $d_j$ is the degree of $F_j$.
Considering  $(x_0,\dots, x_n, t)$ as affine coordinates on
$\C^{n+1}\times \C$ we can read the equations in a different way.
The equations $f_j(x_i)=0$ define the affine cone $C(V)$ over $V$ and
$F_j(x_i,t)=0$ describes a 1-parameter deformation of $C(V)$.
The corresponding infinitesimal deformation is $f_j(x_i)\mapsto
f_j'(x_i)$, which is a deformation of weight $-1$. Conversely,
given a 1-parameter deformation $F_j(x_i,t)=0$ of $C(V)$, with $F_j$
homogeneous of degree $d_j$, we get an extension $W$ of $V$.
For most of the cones considered here the only infinitesimal
deformations of negative weight have weight $-1$ and in that case the
versal deformation in negative weight gives a good description of
all possible extensions.

As the number of equations typically is much larger than the
codimension one needs good ways to describe them.
A prime example is a determinantal scheme $X$:
its ideal is generated by the $t\times t $ minors of an $r\times s$ matrix,
which gives a compact description of the equations.
Following Miles Reid we call this a format.
Canonical curves are not themselves determinantal, but they do
lie on scrolls: a $k$-gonal curve lies on a $(k-1)$-dimensional
scroll, which is given by the minors of a $2\times(g-k+1)$ matrix.
For $k=3$ the curve is a divisor on the scroll, given by one bihomogenous
equation,
and for $k=4$ it is a complete intersection, given by two bihomogeneous
equations. In these cases there is a simple procedure
({\sl `rolling factors'\/}) to write out
one resp.~two sets of equations on $\P^{g-1}$ cutting out the curve
on the scroll.  

Powerful methods exist to compute infinitesimal deformations without
using explicit equations. We  used them for the extension
problem for hyperelliptic curves of high
degree \cite{\rfsth} and trigonal canonical curves \cite{\rfds}.
In these papers also several direct computations with the equations occur.
They seem unavoidable for tetragonal curves,
the subject of a preprint by James N.~Brawner \cite{\rfbrap}. 
The results of these
computations do not depend  on the particular way of choosing the
equations cutting out the curve on the scroll.
This observation was the starting point of this paper. 

We distinguish between different types of deformations and extensions.
If only the equations on the scroll are deformed, but not the scroll
itself we speak of {\sl pure rolling factors deformations}. A typical 
extension lies then on the projective cone over the scroll.
Such a cone is a special case of a scroll of one dimension higher.
If the extension lies on a scroll which is not a cone,
the equations of the scroll are also deformed. We have
a {\sl rolling factors deformation}. Finally if the
extension does not lie on a scroll of one dimension higher
we are in the situation of  a {\sl non-scrollar\/} deformation.
Non-scrollar extensions of tetragonal curves occur only in
connection with Del Pezzo surfaces.
Not every infinitesimal deformation of a scroll gives rise to a deformation
of complete intersections on it. One needs certain lifting 
conditions, which are linear equations in the deformation variables
of the scroll.
Our first main result describes them, 
depending only on the coefficients of the
equations on the scroll.

The next problem is to extend the infinitesimal deformations to a
versal deformation. Here we restrict ourselves to the case that 
all defining equations are quadratic. Our methods thus do not apply 
to trigonal curves, but we can handle tetragonal curves.
Rolling factors obstructions arise.
Previously we observed that one can write them down,
given explicit equations on $\P^n$ \cite{\rfsth, Prop.~2.12}.
Here we give formulas depending only on the coefficients of the
equations on the scroll.
As first application we study base spaces
for hyperelliptic cones. The equations have enough structure
so that explicit solutions can be given.

Surfaces with canonical hyperplane sections are a classical subject.
References to the older literature can be found in Epema's thesis
\cite{\rfep},
which is especially relevant for our purposes. His results say that
apart from $K3$ surfaces only rational surfaces or birationally ruled
surfaces can occur. Furthermore he describes a construction of such
surfaces. Extensions of pure rolling factors type 
of tetragonal curves fit very
well in this description.
A general rolling factors extension is a complete intersection
on a nonsingular four-dimensional scroll. The classification 
of such surfaces \cite{\rfbra}, which we recall below,
shows that surfaces with isolated singularities and in particular
$K3$s can only occur if the degrees of the equations on the scroll
differ at most by 4. 
A tetragonal curve of high genus  with general discrete invariants has
no pure rolling factors deformations. Extensions exist if the
base equations have a solution. 
For low genus we have more variables than equations. For the maximal genus
where almost all curves have a $K3$ extension we find:
 
\proclaim Proposition. The general tetragonal curve
of genus $15$ is hyperplane section of 
$256$ different $K3$ surfaces.

We also look at examples with genus $16$ and $17$.
It is unclear to us which property of a curve makes it have
an extension (apart from the property of being a hyperplane section).

The contents of this paper is as follows.
First we describe the rolling factors format and explain in detail
the equations and relations for the complete intersection
of two divisors on a scroll. Next we recall how canonical curves fit 
into this pattern. In particular we describe the discrete invariants 
for tetragonal curves. The same is done for $K3$ surfaces.
The second section is devoted to the computation of infinitesimal 
deformations. First non-scrollar deformations are treated,
followed by rolling factors deformations. The main result here 
describes the lifting matrix. As application the dimension of $T^1$
is determined for tetragonal cones.
In the third section the base equations for complete intersections
of quadrics on scrolls are derived. As examples base spaces
for hyperelliptic cones are studied. 
The final section describes extensions of tetragonal curves.

\section Rolling factors format.

A subvariety of a determinantal variety can be described by
the determinantal equations and additional equations obtained
by `rolling factors' \cite\rfr. A typical example is the case of 
divisors on scrolls.

We start with a 
$k$-dimensional rational normal scroll $S\subset\P^n$
(for the theory of scrolls we refer to \cite{\rfrch}). The
classical construction is to take $k$ complementary
linear subspaces $L_i$ spanning $\P^n$, each containing  a
parametrised rational normal curve $\phi_i\colon\P^1 \to
C_i\subset L_i$ of degree $d_i=\dim L_i$, and to take
for each $p\in\P^1$ the span of the points $\phi_i(p)$. The
degree of $S$ is $d=\sum d_i=n-k+1$.
If all $d_i>0$ the
scroll $S$ is a $\P^{d-1}$-bundle over $\P^1$.
We allow however that $d_i=0$ for some $i$. Then  
$S$ is the image of $\P^{d-1}$-bundle $\wt S$ over $\P^1$
and $\wt S \to S$ is a rational resolution of singularities.

To give a coordinate description, we take homogeneous
coordinates $(s\cn t)$ on $\P^1$, and 
$(z^{\kl(1)}\cn \cdots\cn z^{\kl(k)})$ on the
fibres. Coordinates on $\P^n$ are $z^{\kl(i)}_{j}=z^{\kl(i)}s^{d_i-j}t^j$,
with $0\leq j\leq d_i$, $1\leq i \leq k$. We give the variable
$z^{\kl(i)}$ the weight $-d_i$.
The scroll $S$ is  given by
the minors of the matrix 
$$
\def\zz#1#2#3{z^{\kl(#1)}_{#2}&\ldots&z^{\kl(#1)}_{#3}}
\Phi=\pmatrix{
\zz10{d_1-1}&\ldots&\zz k0{d_k-1}\cr
\zz11{d_1}&\ldots&\zz k1{d_k}\cr}.
$$

We now consider a divisor on $\wt S$
in the linear system
$|aH-bR|$, where the hyperplane class $H$ and the ruling $R$
generate the Picard group of $\wt S$. When we speak of degree on
$\wt S$ this will be with respect to $H$.
The divisor can be given by 
one bihomogeneous equation $P(s,t,
z^{\kl(i)})$ of degree $a$ in the $z^\kl(i)$, and total degree
$-b$.
By multiplying $P(s,t,z^{\kl(i)})$
with a polynomial of degree $b$ in $(s\cn t)$
we obtain an equation of degree $0$, which can be expressed 
as polynomial of degree $a$ in the $z^{\kl(i)}_j$; this expression
is not unique, but the difference of two expressions lies in the
ideal of the scroll. By the obvious choice,
multiplying with $s^{b-m}t^m$, we obtain $b+1$ equations
$P_m$. In the transition
from the equation $P_m$ to $P_{m+1}$ we have to increase  by
one the sum of the lower indices of the factors $z^{\kl(i)}_j$ in each
monomial, and we can and will always achieve this by increasing 
exactly one index. This amounts to replacing a $z^{\kl(i)}_j$, which
occurs in the top row of the matrix, by the element
$z^{\kl(i)}_{j+1}$ in the bottom row of the same column. This is the
procedure of `{\sl rolling factors\/}'.

\roep Example \num.
Consider the cone over $2d-b$ points in $\P^d$,
lying on a rational normal curve of degree $d$, with $b<d$.
Let the polynomial $P(s,t)= p_0s^{2d-b}+p_1s^{2d-b-1}t+
\cdots+p_{2d-b}t^{2d-b}$
determine the points on the rational curve.   
We get the determinantal
$$
\bmatrix{
z_0&z_1&\ldots&z_{d-1}\cr
z_1&z_2&\ldots&z_{d}\cr} 
$$ 
and additional equations $P_m$. To be specific we assume that $b=2c$:
$$
\matrix{
P_0&=&
       \hfill p_0z_0^2+p_1z_0z_1+{}&\hskip-1em\cdots&\hskip-1em{}+
         p_{2d-2c-1}z_{d-c-1}z_{d-c}+p_{2d-2c}z_{d-c}^2\hfill\cr
P_1&=&
      \hfill   p_0z_0z_1+p_1z_1^2+{}&\hskip-1em\cdots&\hskip-1em{}+
           p_{2d-2c-1}z_{d-c}^2+p_{2d-2c}z_{d-c}z_{d-c+1}\hfill\cr
&\vdots&\cr
P_{2c}&=&
       \hfill   p_0z_c^2+p_1z_cz_{c+1}+{}&\hskip-1em\cdots&
            \hskip-1em{}+p_{2d-2c-1}z_{d-1}z_{d}+p_{2d-2c}z_{d}^2\;.\hfill \cr
} 
$$
\endroep

The `rolling factors' phenomenon can also occur  if the entries
of the matrix are more general.

\roep Example \num.
Consider a non-singular hyperelliptic curve of genus 5, with
a half-canonical line bundle $L=g_1^2+P_1+P_2$ where the $P_i$ 
are Weierstrass
points. According to \cite{\rfr}, Thm.~3,
the ring $R(C,L)=\bigoplus H^0(C,nL)$  is 
$k[x_1,x_2,y_1,y_2,z_1,z_2]/I$ with $I$ given by the determinantal
$$
\bmatrix{x_1&y_1&x_2^2&z_1\cr x_2& x_1^2&y_2&z_2}
$$
and the three rolling factors equations
$$
\eqalign{
z_1^2  &=  x_1^2h  +y_1^3 +x_2^4 y_2  \cr
z_1z_2  &=  x_1x_2h  +y_1^2x_1^2 +x_2^2 y_2^2  \cr
z_2^2  &=  x_2^2h  +y_1x_1^4 + y_2^3  \cr
}
$$
where $h$ is some quartic in $x_1$, $x_2$, $y_1$, $y_2$.
\endroep

The description of the syzygies of  a subvariety
$V$ of the scroll $S$ proceeds in two steps. First one constructs
a resolution of $\sier {\wt V}$ by vector bundles on $\wt S$
which are repeated extensions of line bundles.
Schreyer describes,
following Eisenbud,  Eagon-Northcott type complexes ${\cal
C}^b$ such that ${\cal C}^b(a)$ is the minimal resolution
of $i_*\big(\sier {\wt S}(-aH+bR)\big)$ as $\sier{\P^n}$-module, if $b\geq-1$
\cite{\rfschr}. Here $i\colon \wt S \to \P^n$ is the map defined by $H$.
The resolution of $\sier V$ is then
obtained by taking an (iterated)  mapping cone.

The matrix $\Phi$ defining the scroll
can be obtained intrinsically from the multiplication map 
$$
H^0\sier{\wt S}(R)\otimes H^0\sier{\wt S}(H-R) \lra H^0\sier{\wt S}(H)\;.
$$
In general, given a map $\Phi\colon F \to G$ of locally free sheaves of rank
$f$ and $g$ respectively, $f\geq g$, on a variety one defines 
Eagon-Northcott type complexes ${\cal C}^b$, $b\geq-1$,
in the following way:
$$
{\cal C}_j^b=\cases{\bigwedge^j F\otimes S_{b-j}G, &for $0\leq j \leq b$\cr
\bigwedge^{j+g-1} F\otimes D_{j-b-1}G^*\otimes \bigwedge^gG^*,
&for $j\geq b+1$\cr}
$$
with differential defined by multiplication with $\Phi\in F^*\otimes G$
for $j\neq b+1$ and $\bigwedge^g\Phi \in \bigwedge^g F^* 
\otimes \bigwedge^g G$ for $j=b+1$ in the appropriate term
of the exterior, symmetric or divided power algebra.

In our situation $F\cong \sier{\P^n}^d(-1)$ and $G\cong \sier {\P^n}^2$ with
$\Phi$ given by the matrix of the scroll. Then ${\cal C}^b(-a)$
is for $b\geq -1$ the minimal resolution of 
$\sier {\wt S}(-aH+bR)$ as $\sier\P$-module \cite{\rfschr, Cor.~1.2}.

Now let $V\subset S \subset \P^n$ be a `complete intersection'
of divisors $Y_i\sim a_iH -b_iR$, $i=1,$ \dots, $l$, on 
a $k$-dimensional rational scroll of degree $d$ with $b_i\geq 0$.
The resolution of $\sier V$ as $\sier S$-module is a
Koszul complex and the iterated mapping cone of complexes ${\cal C}^b$
is the minimal resolution \cite{\rfschr, Sect.~3, Example}.

To make this resolution more explicit we look at the case 
$l=2$, which is  relevant for tetragonal curves. The iterated
mapping cone is
$$
\big[ \cc^{b_1+b_2}(-a_1-a_2) \lra 
   \cc^{b_1}(-a_1) \oplus \cc^{b_2}(-a_2) \big]
 \lra \cc^0
$$
To describe equations and relations we give the first steps
of this complex. 
We first consider the case that $b_1\geq b_2 >0$.
We write $\so$ for $\sier {\P^n}$. We get the double complex
$$

\matrix{\so&\leftarrow &\bigwedge^2\so^d(-1)\quad
\leftarrow \quad\bigwedge^3\so^d(-1)\otimes\so^2 \cr
\Big\uparrow & &\Big\uparrow \hskip 4cm\cr
S_{b_1}\so^2(-a_1)\oplus S_{b_2}\so^2(-a_2)&\leftarrow&
\so^d(-1)\otimes S_{b_1-1}\so^2(-a_1)
  \oplus \so^d(-1)\otimes S_{b_2-1}\so^2(-a_2)
\cr\Big\uparrow \cr
S_{b_1+b_2}\so^2(-a_1-a_2) \cr
}
$$
The equations for $V$ consist of the determinantal ones plus
two sets of additional equations obtained by rolling factors:
the two equations $P^\kl(1)$, $P^\kl(2)$ defining $V$ on the scroll 
give rise to $b_1+1$ equations $P^\kl(1)_m$
and $b_2+1$ equations $P^\kl(2)_m$.

To describe the relations 
we introduce the following notation.
A column in the matrix $\Phi$ has the form $(z^{\kl(i)}_j, 
z^{\kl(i)}_{j+1})$. We write symbolically
$(z_\al,z_{\al+1})$, where the index $\al$ stands for the pair
${}^{\kl(i)}_{j}$  and $\al+1$ means adding $1$ to the lower index.
More generally, if $\al={}^{\kl(i)}_{j}$ and  $\al'={}^{\kl(i')}_{j'}$
then the sum $\al+\al':=j+j'$ only involves the lower indices.
To access the upper index we say that $\al$ is of type $i$.
The rolling factors assumption is that 
two consecutive
additional equations are of the form
$$
\eqalign{
P_m=&\sum_\al p_{\al,m}z_\al,\cr
P_{m+1}=&\sum_\al p_{\al,m}z_{\al+1}.\cr}
$$
where the polynomials $p_{\al,m}$ depend on the $z$-variables
and the sum runs over all possible pairs $\al={}^{\kl(i)}_{j}$.
To roll from $P_{m+1}$ to $P_{m+2}$ we collect the `coefficients'
in the equation $P_{m+1}$ in a different way: we also have
$P_{m+1}=\sum_\al p_{\al,m+1}z_\al$.

We write the scrollar equations as 
$f_{\al\be}=z_\al z_{\be+1}-z_{\al+1}z_\be$.
The relations between them  are
$$
\displaylines{
R_{\al,\be,\ga}=f_{\al,\be}z_\al-f_{\al,\ga}z_\be+f_{\be,\ga}z_\al,\cr
S_{\al,\be,\ga}=f_{\al,\be}z_{\al+1}-f_{\al,\ga}z_{\be+1}
                    +f_{\be,\ga}z_{\al+1},\cr}
$$
which corresponds to the term 
$\bigwedge^3\so ^d(-1) \otimes \so ^2$ in Schreyer's resolution.
The second line yields
relations involving the two sets of $P^\kl(n)_m\,$:
$$
R^n_{\be,m}=P^\kl(n)_{m+1}z_\be-P^\kl(n)_m z_{\be+1}-
\sum_\al  f_{\be,\al}p^\kl(n)_{\al,m},
$$
where $n=1,2$ and $0\leq m <b_i$. 
We note the following relation: 
$$
R^n_{\be,m}z_\ga-R^n_{\be,m}z_\be-
\sum R^n_{\be,\ga,\al}p^\kl(n)_{\al,m}=
P^\kl(n)_m f_{\be,\ga}-f_{\be,\ga}P^\kl(n)_m.
$$
The right hand side is a Koszul relation; the second factor in
each product is considered as coefficient. There are similar
expressions involving $z_{\ga+1}$, $z_{\be+1}$ and $S_{\be,\ga,\al}$.
Finally by multiplication with suitable powers
of $s$ and $t$ 
the Koszul relation $P^\kl(1)P^\kl(2)-P^\kl(2)P^\kl(1)$ gives
rise to $b_1+b_2+1$ relations --- this is the term
$S_{b_1+b_2}\so^2(-a_1-a_2)$.

In case $b_1>b_2=0$ the resolution is
$$

\matrix{\so&\leftarrow &\bigwedge^2\so^d(-1)\quad
\leftarrow \quad\bigwedge^3\so^d(-1)\otimes\so^2 \cr
\Big\uparrow & &\Big\uparrow \hskip 4cm\cr
S_{b_1}\so^2(-a_1)\oplus \so^2(-a_2)&\leftarrow&
\so^d(-1)\otimes S_{b_1-1}\so^2(-a_1)
  \oplus \bigwedge^2\so^d(-1-a_2)
\cr\Big\uparrow \cr
S_{b_1}\so^2(-a_1-a_2) \cr
}
$$
The new term expresses the Koszul relations between the one equation
$P^\kl(2)$ and the determinantal equations 
(which had previously been expressible in terms of rolling factors
relations). For the computation of deformations these relations
may be ignored.

Finally, if $b_2=-1$, the equations change drastically.

\roep {\rm(\num)} Canonical curves {\rm\cite\rfschr}.

\noindent 
A $k$-gonal canonical curve lies on a $(k-1)$-dimensional scroll
of degree $d=g-k+1$. We write $D$ for the divisor of the $g^1_k$.
To describe the type $S(e_1,...,e_{k-1})$ of the scroll we introduce 
the numbers
$$
f_i=h^0(C,K-iD)-h^0(C,K-(i+i)D) = k+h^0(iD)-h^0((i+1)D)
$$
for $i\geq0$ and set 
$$
e_i=\#\{j\mid f_j\geq i\}-1\;.
$$
In particular, $e_1$ is the minimal number $i$ such that
$h^0((i+1)D)-h^0(iD)=k$ and it satisfies therefore $e_1\leq \frac{2g-2}k$.

A trigonal curve lies on a scroll 
of type $S(e_1,e_2)$ and degree $d=e_1+e_2=g-2$ with
$$
\frac{2g-2}3 \geq e_1 \geq e_2 \geq \frac{g-4}3
$$
as a divisor of type $3H-(g-4)R$.
The minimal resolution of $\sier C$ is given by the mapping cone
$$
\cc^{d-2}(-3)\lra \cc^0\;.
$$
\global\edef\labsect{{\the\secno.\the\thmno}}
Introducing bihomogeneous coordinates $(x\cn y;s\cn t)$ and coordinates
$x_i=xs^{e_1-i}t^i$, $y_i=ys^{e_2-i}t^i$ we obtain the scroll
$$
\pmatrix{
x_0&x_1&\ldots&x_{e_1-1}&y_0&y_1&\ldots&y_{e_2-1}\cr
x_1&x_2&\ldots&x_{e_1  }&y_1&y_2&\ldots&y_{e_2}
}
$$
and a bihomogeneous equation for $C$
$$
P=A_{2e_1-e_2+2}x^3+B_{e_1+2}x^2y+C_{e_2+2}xy^2+D_{2e_2-e_1+2}y^3
$$
where $A_{2e_1-e_2+2}$ is a polynomial in $(s\cn t)$ of degree  
$2e_1-e_2+2$ and similarly for the other coefficients.
By rolling factors $P$ gives rise to $g-3$ extra equations.

The inequality $e_1\leq \frac{2g-2}3$ can also be explained from
the condition that the curve $C$ is nonsingular, 
which implies that the polynomial $P$ is irreducible,
and therefore  the degree $2e_2-e_1+2=2g-2-3e_1$
of the polynomial $D_{2e_2-e_1+2}$ is nonnegative. 
The other inequality follows from this one
because $e_1=g-e_2-2$, but also by considering the degree of $A_{2e_1-e_2+2}$.

A tetragonal curve of genus $g\geq 5$ is a complete intersection of
divisors $Y\sim 2H-b_1R$ and $Z\sim 2H-b_2R$ on a scroll of
type $S(e_1,e_2,e_3)$ of degree $d=e_1+e_2+e_3=g-3$, with $b_1+b_2=d-2$,
and
$$
\frac {g-1}2\geq e_1 \geq e_2 \geq e_3 \geq 0
$$
We introduce bihomogeneous coordinates $(x\colon y\cn z;s\cn t)$.
Then $Y$ is given by an equation
$$
P=P_{1,1}x^2+P_{1,2}xy+\cdots+P_{3,3}z^2
$$
with $P_{ij}$ (if nonzero) a polynomial in $(s\cn t)$ of degree $e_i+e_j-b_1$ 
and likewise $Z$ has equation 
$$
Q=Q_{1,1}x^2+Q_{1,2}xy+\cdots+Q_{3,3}z^2
$$
with $\deg Q_{ij}=e_i+e_j-b_2$.

The minimal resolution is of type discussed above, because
the condition $-1\leq b_2\leq b_1 \leq d-1$ is satisfied: the only 
possibility   to have a divisor of type $2H-bR$ with $b\geq d$ is to
have $e_1=e_2=d/2$, $e_3=0$ and $b=d$, but then the equation $P$ is 
of the form $\al x^2+\be xy +\ga y^2$ with constant coefficients, so
reducible.
If $b_2=-1$ also cubics are needed to generate the ideal, so the 
curve admits also a $g_3^1$ or $g_5^2$; this happens only up to
$g=6$. We exclude these cases and assume that $b_2\geq 0$.

\proclaim Lemma \num.
We have $b_1\leq 2e_2$ and $b_2\leq 2e_3$.

\roep Proof.
If $b_1> 2e_2$ the polynomials $P_{22}$, $P_{23}$ and $P_{33}$ vanish
so $P$ is reducible and therefore $C$.
If $b_2> 2e_3$ then $P_{33}$ and $Q_{33}$ vanish. This means that the
section $x=y=0$ is a component of $Y\cap Z$ on the $\P^2$-bundle
whose image in $\P^{g-1}$ is the scroll (if $e_3>0$ the scroll is nonsingular,
but for $e_3=0$ it is a cone). As the arithmetic genus of $Y\cap Z$ is
$g$ and its image has to be the nonsingular curve $C$ of genus $g$,
the line cannot be a component.
\qed

This Lemma is parts 2 -- 4 in \cite{\rfbra, Prop.~3.1}.
Its last part is incorrect. It states that $b_1 \leq e_1+e_3$ if $e_3>0$,
and builds upon the fact that $Y$ has only isolated singularities. 
However the discussion in \cite{\rfschr} makes clear that this need not be the
case.

The surface $Y$ fibres over $\P^1$. There are now two cases, first that
the general fibre is a non-singular conic. In this case one of the 
coefficients $P_{13}$, $P_{23}$ or $P_{33}$ is nonzero, giving indeed
$b_1\leq e_1+e_3$. 

The other possibility is that each fibre is a singular conic. Then
$Y$ is a birationally ruled surface over a (hyper)elliptic curve
$E$ with a rational curve $\wl E$ of double points, the canonical
image of $E$, and $C$ does not intersect $\wl E$. 
This means that the section $\wl E$ of the scroll
does not intersect the surface $Z$, so if one inserts the parametrisation
of $\wl E$ in the equation of $Z$ one obtains a non-zero constant.
Let the section be given by polynomials in $(s\cn t)$, which
if nonzero have degree  $d_s-e_1$, $d_s-e_2$, $d_s-e_3$. Inserting them
in the polynomial $Q$ gives a polynomial of degree $2d_s-b_2$.
So $b_2$ is even and $2d_s=b_2\leq 2e_3$. On the other hand
$d_s-e_3\geq 0$ so $d_s=e_3$ and $b_2=2e_3$.
The genus of $E$ satisfies $p_a(E)=b_2/2+1$.
If $b_1>e_1+e_3$, then $Y$ is singular along the section $x=y=0$.
An hyperelliptic involution
can also occur if $b_1\leq e_1+e_3$.

We have shown:
\proclaim Lemma \num.
If $Y$ is singular, in particular if $b_1>e_1+e_3$, then $b_2= 2e_3$.

Finally we analyse the case $b_2=0$ (cf.~\cite{\rfbrap}).

\proclaim Lemma \num.
A nonsingular tetragonal curve is bielliptic or lies on a Del Pezzo surface
if and only if $b_2=0$. The first case occurs for $e_3=0$, and the second
for the values $(2,0,0)$, $(1,1,0)$, $(2,1,0)$, $(1,1,1)$, $(3,1,0)$,
$(2,2,0)$, $(2,1,1)$, $(3,2,0)$, $(2,2,1)$, $(4,2,0)$, $(3,2,1)$ or $(2,2,2)$
of the triple $(e_1,e_2,e_3)$.

\roep Proof.
If the curve is bielliptic or lies on a Del Pezzo, the $g_4^1$ is not unique,
which implies that the scroll is not unique. This is only possible if $b_2=0$
by \cite{\rfschr}, p.~127. Then $C$ is the complete intersection of a quadric
and a surface $Y$ of degree $g-1$, which is uniquely determined
by $C$. 

The inequality $e_1+e_2+e_3-2=b_1\leq 2e_2$ shows that 
$e_3\leq e_2-e-1+2\leq 2$. 
If the general fibre of $Y$ over $\P^1$ is non-singular
we have $b_1\leq e_1+e_3$.
This gives $e_2\leq 2$ and $b_1\leq 4$.
The possible values are now easily determined.
If the general fibre of $Y$ is singular then  $e_3=b_2/2=0$ and 
$Y$ is an elliptic cone.
\qed

\roep {\rm (\num)} $K3$ surfaces.

\noindent
Let $X$ be a $K3$ surface (with at most rational double point
singularities)  on a scroll. If the scroll is 
nonsingular the projection onto $\P^1$ gives an elliptic
fibration on $X$, whose general fibre is smooth.
This is even true if the scroll is singular: the strict transform
$\wt X$ on $\wt S$ has only isolated singularities.

We start with the case of divisors. 
A treatment of such 
scrollar surfaces with an elliptic fibration
can be found in \cite{\rfrch, 2.11}. One finds:

\proclaim Lemma \num.
For the general $F\in |3H-kR|$ on a scroll $S(e_1,e_2,e_3)$
the general fibre of the
elliptic fibration is a nonsingular  cubic curve
if and only  if
$k\leq 3e_2$ and $k \leq e_1+2e_3$. 

\comment
\roep Proof.
We take bihomogeneous coordinates
$(x\colon y\cn z;s\cn t)$ in which  $F$ is given by a bihomogeneous 
equation
$$
A_{3e_1-k}(s,t)x^3+ \cdots + A_{3e_3-k}(s,t)z^3
$$  
The condition that this cubic is irreducible means in particular that the
equation is not divisible by $x$, i.e. at least the term $y^3$ occurs and
the degree of $A_{3e_2-k}$ is nonnegative.
The condition that $(0\cn 0\cn 1)$ is not  singular point is one on
the coefficient of $xz^2$: it gives $\deg A_{e_1+2e_3-k}\geq0$. 
As the curve $x^3+y^3+xz^2$ is smooth, these conditions are
also sufficient for non-singularity of the general fibre.
\qed
\endcomment

If one fixes $k$ and $e_1+e_2+e_3$ these conditions limit the possible
distribution of the integers $(e_1,e_2,e_3)$.
By the adjunction formula one has $k=e_1+e_2+e_3-2$ for a $K3$ surface.
In this case we
obtain 12 solutions, which fall into 3 deformation types of scrolls, 
according to
$\sum e_i \pmod 3$:
$$
\displaylines{
(e+2,e,e-2) \to ( e+1,e,e-1) \to (e,e,e)  \cr
(e+3,e,e-2)\to(e+2,e,e-1)\to(e+1,e+1,e-1)\to (e+1,e,e)\cr
(e+4,e,e-2) \to  (e+3,e,e-1) \to (e+2,e+1,e-1) \to (e+2,e,e) \to (e+1,e+1,e)
\cr
}
$$ 
The general element of the linear system can only have singularities at the
base locus. The base locus is the section $(0\cn 0\cn 1)$ if and only if
$k>3e_3$ and there
is a singularity at the points $(s\cn t)$ where both $A_{e_1+2e_3-k}$ and 
$A_{e_2+2e_3-k}$ vanish. The assumption that the coefficients
are general implies now that 
$\deg A_{e_2+2e_3-k}<0$ and $\deg A_{e_1+2e_3-k}>0$.

In the 12 cases above this occurs only for $(e+3,e,e-1)$ and 
$(e+2,e,e-1)$. In the first case the term $y^2z$ is also missing, yielding
that there is an $A_2$-singularity at the only zero of $A_{e_2+2e_3-k}$,
whereas the second case gives an $A_1$.
The scroll $S_{e+4,e,e-2}$ deforms into $S_{e+3,e,e-1}$, but the general 
$K3$-surface on it does not deform to a $K3$ on $S_{e+3,e,e-1}$, but only
those with an $A_2$-singularity.
These results hold if all $e_i>0$; we leave the modifications
in case $e_3=0$ to the reader.

\comment
These results hold if all $e_i>0$, and have to be modified if $e_3=0$.
We consider the scroll 
$S_{e+4,e,e-2}$ with $e=2$, i.e. $S_{6,2,0}$.
The general element of the linear system  still is
$$
\vp_{12}(s,t)x^3+\vp_8(s,t)x^2y+\vp_6(s,t)x^2z+\vp_4(s,t)xy^2+\vp_2(s,t)xyz
+\vp_0xz^2+\psi_0y^3
$$
with $\vp_i$ homogeneous polynomials of degree $i$ in $(s\cn t)$ and $t$. 
Along the section $(0\cn 0\cn 1)$  of the $\P^2$-bundle $\wt S$
the linear part of the equation is 
$\vp_0x=0$ and therefore the normal bundle of the curve in the
resolution of the $K3$-surface is the same as 
of the curve in the resolution of the two-dimensional scroll $S_{2,0}$,
showing that the curve is a $(-2)$-curve. The singular
$K3$-surface has therefore an ordinary double point. The same holds for 
$S_{e+3,e,e-2}$ and $S_{e+2,e,e-2}$ with $e=2$, and
$S_{e+2,e+1,e-1}$ and $S_{e+1,e+1,e-1}$ with $e=1$.
For $S_{e+3,e,e-1}$ and $S_{e+2,e,e-1}$ with $e=1$ we have as
linear part $\vp_1x=0$: the normal bundle of the section $(0\cn 1)$ 
in $\wt S_{1,0}$ has degree $-1$, but the zero of $\vp_1$ (at the singularity
of the singular $K3$) makes the curve a $(-2)$-curve on the resolution of
the $K3$. Finally for $\wt S_{2,1,0}$ we have linear part $\vp_1x+\vp_0y$,
making the curve again a $(-2)$.
In this case we make the rolling factors equations explicit.
We have in the $\P^5$ with coordinates 
$(x_0\cn x_1\cn x_2\cn y_0\cn y_1\cn z)$ the scroll
$$
\pmatrix{x_0&x_1&y_0 \cr x_1&x_2&y_1 }
$$ 
and equations 
$$
\displaylines{
y_0z^2+ p_1 y_0^2z + p_2 y_0y_1z + \cdots + p_lx_1x_2^2 \cr 
y_1z^2+ p_1 y_0y_1z +p_2 y_1^2z + \cdots + p_lx_2^3 \;.\cr}
$$
Using inhomogeneous coordinates with $z=1$ we can
eliminate $y_0$ and $y_1$ by means of the two additional equations,
showing that we indeed have a singularity of type $x_0x_2-x_1^2$.
\endcomment

\pageinsert
$$
\vbox{
\def\twee{height1.5pt&\omit&&\omit&&\omit&&\omit&&\omit&\cr}
\offinterlineskip
\hrule
\halign{&\vrule#&\strut\quad  $\hfil #\hfil$ \quad\cr
\twee
&(b_1,b_2)&& (e_1,e_2,e_3,e_4)&& \hbox{\# mod.}&& \hbox{base} && \hbox{sings} &\cr
\twee
\noalign{\hrule}
\twee
\noalign{\hrule}
\twee
&(2e,2e-2)   && (e+3,e+1,e-1,e-3)  && 17 &&  B_{e-1} &&      -- & \cr  
&            && (e+3,e,e-1,e-2) && 15 &&    B_{e-1} &&     A_3  & \cr
 &           && (e+2,e+1,e-1,e-2) && 16 &&   B_{e-1} &&     A_1  & \cr   
&            && (e+2,e,e,e-2) && 16 &&    B_{e-2} &&           -- & \cr 
&            && (e+2,e,e-1,e-1) && 15 &&   B_{e-1} &&     2A_1  & \cr
&           && (e+1,e+1,e-1,e-1) && 16 &&    B_{e-1} &&         -- & \cr
&           && (e+1,e,e,e-1) && 17 &&    B_{e-1} &&          -- & \cr
&           && (e,e,e,e) && 17 &&  \emptyset &&          --  & \cr
\twee
\noalign{\hrule}
\twee
&(2e-1,2e-1)  && (e+1,e+1,e,e-2) && 17 &&   B_{e-2} && --& \cr
&             && (e+1,e,e,e-1) && 17 &&    B_{e-1} && --& \cr
&             && (e,e,e,e) && 18 &&  \emptyset &&        --& \cr
\twee
\noalign{\hrule}
\twee
\noalign{\hrule}
\twee
&(2e+1,2e-2)   && (e+4,e+1,e-1,e-3) && 17 &&   B_{e-1} &&  -- & \cr
&            && (e+3,e+1,e-1,e-2) && 16 &&    B_{e-1} &&   A_1 & \cr
&            && (e+2,e+1,e-1,e-1) && 16 &&    B_{e-1} &&    --& \cr
&             && (e+1,e+1,e,e-1) && 17 &&    B_{e} &&    --& \cr
\twee
\noalign{\hrule}
\twee
&(2e,2e-1)   && (e+2,e+1,e,e-2) && 17 &&    B_{e-2} &&  --& \cr
&           && (e+2,e,e,e-1) && 15 &&      B_{e-1} &&  A_1 & \cr
&           && (e+1,e+1,e,e-1) && 17 &&    B_{e-1} &&  --& \cr
&           && (e+1,e,e,e) && 18 &&    \emptyset &&        --& \cr
\twee
\noalign{\hrule}
\twee
\noalign{\hrule}
\twee
&(2e+2,2e-2)  && (e+5,e+1,e-1,e-3) && 18 &&    B_{e-1} &&   -- & \cr
&            && (e+4,e+1,e-1,e-2) && 17 &&    B_{e-1} &&   A_1    & \cr
&             && (e+3,e+1,e-1,e-1) && 17 &&    B_{e-1} &&      --& \cr
&             && (e+2,e+1,e,e-1) && 18 &&    B_{e} &&    --& \cr
&            && (e+1,e+1,e+1,e-1) && 18 &&    B_{e-1} &&    --& \cr
\twee
\noalign{\hrule}
\twee
&(2e+1,2e-1)  && (e+3,e+1,e,e-2) && 16 &&    B_{e} &&   -- & \cr
&           && (e+2,e+1,e,e-1) && 16 &&    B_{e} &&     --  & \cr
&            && (e+1,e+1,e,e) && 17 &&     B_{e} &&   --& \cr
\twee
\noalign{\hrule}
\twee
&(2e,2e)   && (e+2,e+2,e,e-2) && 17 &&    B_{e-2} &&   -- & \cr
&            && (e+2,e+1,e,e-1) && 16 &&    B_{e-1} &&  A_1 & \cr
&          && (e+1,e+1,e+1,e-1) && 17 &&      B_{e-1} &&   -- & \cr
&            && (e+2,e,e,e) && 15 &&    \emptyset &&        -- & \cr
&         && (e+1,e+1,e,e) && 17 &&     \emptyset &&        --& \cr
\twee
\noalign{\hrule}
\twee
\noalign{\hrule}
\twee
&(2e+2,2e-1)   && (e+4,e+1,e,e-2) && 16 &&    B_{e} && A_1  & \cr
&          && (e+3,e+1,e,e-1) && 16 &&    B_{e} &&    A_1  & \cr
&            && (e+2,e+1,e,e) && 17 &&    B_{e} &&  A_1  & \cr
&            && (e+1,e+1,e+1,e) && 17 &&     B_{e} &&   A_1  & \cr
\twee
\noalign{\hrule}
\twee
&(2e+1,2e)   && (e+3,e+2,e,e-2) && 17 &&    B_{e} &&     -- & \cr
&          && (e+3,e+1,e,e-1) && 15 &&   B_{e} &&   A_2  & \cr
&             && (e+2,e+2,e,e-1) && 16 &&   B_{e} &&  A_1 & \cr
&             && (e+2,e+1,e+1,e-1) && 17 &&   B_{e-1} && --& \cr
&            && (e+2,e+1,e,e) && 16 &&     B_{e} &&    --& \cr
&            && (e+1,e+1,e+1,e) && 18 &&     B_{e} &&    --& \cr
\twee}
\hrule
}
$$
\endinsert

The tetragonal case is given as exercise in \cite{\rfrch} 
and the complete solution (modulo some minor mistakes)
can be found in \cite{\rfbra}. We give the results:

\proclaim Lemma \num.
For the general complete intersection of divisors of type $2H-b_1R$ and
$2H-b_2R$ on a scroll $S_{e_1,e_2,e_3,e_4}$ 
the general fibre of the
elliptic fibration is a nonsingular  quartic curve
if and only if either
\item{$\al$:} $b_1\leq e_1+e_3$, $b_1\leq 2e_2$ and $b_2\leq 2e_4$, or
\item{$\be$:}  $b_1\leq e_1+e_4$, $b_1\leq 2e_2$, $2e_4<b_2\leq 2e_3$
and $b_2\leq e_2+e_4$.  

\comment
\roep Proof.
Let the surface be the intersection of divisors $Y$ and $Z$.
We take bihomogeneous coordinates $(x\cn y\cn z\cn w;s\cn t)$.
The condition that the fibre of $Y$ is irreducible 
implies firstly that the equation
is not divisible by $x$, so the coefficient of $y^2$ has nonnegative
degree: $b_1\leq 2e_2$, and furthermore that  the equation does not only
involve $x$ and $y$, giving: $b_1\leq e_1+e_3$. 
That the fibre of $Y\cap Z$ is irreducible implies that the curve $x=y=0$ is  
not a component,  giving $b_2\leq 2e_3$, and also that it is not a cone:
$b_2\leq e_1+e_4$.
Finally, the condition that $(0\cn 0\cn 0\cn 1)$ is not a singular point
can be satisfied by requiring that this point does not lie on the 
curve: $b_2\leq 2e_4$, or by requiring that the inhomogeneous linear parts
are linearly independent: the monomial $xw$ occurs for $Y$ ($b_1\leq e_1+e_4$)
and $yw$ for $Z$ ($b_2\leq e_2+e_4$). By removing redundant conditions
we arrive at necessary conditions. It is easy to construct smooth curves
with the necessary monomials, showing also sufficiency.
\qed

We determine the singularities of the general element. By Bertini's
theorem singularities can only occur at points in the base locus of
the linear systems $|2H-b_1R|$ and
$|2H-b_2R|$.
The base locus is a subscroll, for which we use the following notation
\cite{\rfrch}, 2.8.
For any  $a$ we denote by $B_a$
the subscroll corresponding to the subset of all $e_i$ with $e_i\leq a$,
defined by the equations $z^{\kl(j)}=0$ for $e_j>a$.
\endcomment

\def\itemind{\par\hskip -10pt\hang\textindent}
\proclaim Proposition \num.
The general element is singular at a point of 
the section $(0\cn 0\cn 0\cn 1)$ if the invariants satisfy
in addition one of the following conditions:
\itemind{$1\al$:} $b_2<2e_4$, $b_1>e_1+e_4$.
\itemind{$1\be i$:} $b_2\leq e_3+e_4$, $e_2+e_4<b_1<e_1+e_4$.
\itemind{$1\be ii$:}  $b_2> e_3+e_4$, and $e_1+e_2+2e_4>b_1+b_2$.
\hfil\break
There is a singularity with $z\neq0$ if
\itemind{$2\al$:}  $b_1>e_2+e_3$,  $e_1+e_3>b_1>e_1+e_4$.
\itemind{$2\al\be$:} $e_1+e_4\geq b_1>e_2+e_3$ and
\itemitem{i:} if $b_2\leq 2e_4$ then $2(e_1+e_3+e_4)>2b_1+b_2$
\itemitem{ii:} if $2e_4<b_2\leq e_3+e_4$ then $e_1+2e_3+e_4>b_1+b_2$
\itemitem{iii:} $e_3+e_4<b_2<2e_3 $ 

\comment
\roep Proof.
We look at the two bihomogeneous equations and form the Jacobian matrix
with respect to all variables $(x\cn y\cn z\cn w;s\cn t)$.

Suppose there is a singularity at a point of the section $(0\cn 0\cn 0\cn 1)$.
If this section is not the base locus of $|Z|$, 
then the point can only lie
at a zero of the coefficient $Q_{4,4}(s,t)$ of $w^2$, giving
$b_2<2e_4$. The threefold $Z$
is smooth at such points and the  general $Y$ is transverse to it, if
smooth. So the necessary and sufficient condition for a singularity
is that the coefficient of $xw$ in the equation of $Y$ vanishes:
$b_1>e_1+e_4$.

If the section lies in the base locus (i.e.~$b_2>2e_4$), then the condition is that
the linear parts of the equations in the  affine chart $w=1$ are linearly
dependent:
$$
\Rank \pmatrix {P_{1,4}& P_{2,4}&P_{3,4}\cr
               Q_{1,4}& Q_{2,4}&Q_{3,4}\cr}\leq1
$$
The minor $P_{1,4}Q_{2,4}-P_{2,4}Q_{1,4}$ does not vanish identically,
because $P_{1,4}$ and $Q_{2,4}$ don't.
If $P_{2,4}Q_{3,4}-P_{3,4}Q_{2,4}$ does not vanish identically,
both polynomials will in general not have common roots.
If $Q_{3,4}\not\equiv 0$ then $P_{2,4}\equiv 0$ 
(which implies $P_{3,4}\equiv0$) and 
the singularities are the zeroes of $P_{1,4}$.
Otherwise
$Q_{3,4}\equiv 0$ (which also implies $P_{3,4}\equiv0$) 
and the number of singularities
is the degree of $P_{1,4}Q_{2,4}-P_{2,4}Q_{1,4}$.

We now assume that there is a singularity with $z\neq0$.
The base locus of $Y$ is now $x=y=0$, so in the first line of the
Jacobian matrix
the derivatives w.r.t $z$, $w$, $s$ and $t$ vanish identically.
The curve $Q_{3,3}z^2+Q_{3,4}zw+Q_{4,4}w^2=0$
has for general values of the parameters no singularities,
so the partial derivatives have no common zeroes.
The condition for a singularity is therefore that the first
line of the Jacobian matrix vanishes, i.e.,
$P_{1,3}z+P_{1,4}w=P_{2,3}z+P_{2,4}w=0$; that the singular point lies on the
surface $Z$ gives the condition that $Q_{3,3}z^2+Q_{3,4}zw+Q_{4,4}w^2=0$.
If $P_{2,3}\not\equiv 0$ the first two conditions give only finitely
many solutions $(s\cn t;z\cn w)$ with $z\neq0$, which will not satisfy the third
equation. So  $P_{2,3}\equiv 0$. 
Suppose first that $P_{1,4}\equiv0$ (which implies that 
$Q_{4,4}\not\equiv0$). Then there are two singular points for each zero 
$(s_0\cn t_0)$ of $P_{1,3}$, determined by the two pairs $(z\cn w)$ which 
satisfy $Q_{3,3}(s_0,t_0)z^2+Q_{3,4}(s_0,t_0)zw+Q_{4,4}(s_0,t_0)w^2=0$.
Otherwise we use the equation $P_{1,3}z+P_{1,4}w=0$
to eliminate $w$ and  find as condition for a singular point with $z\neq0$
that $Q_{3,3}P_{1,4}^2-Q_{3,4}P_{1,4}P_{1,3}+Q_{4,4}P_{1,3}^2=0$.
If $Q_{4,4}\equiv0$ then a zero of $P_{1,4}$ gives a singularity
at the section $(0\cn 0\cn 0\cn 1)$, so the condition for a singular point
with $z\neq 0$ reduces to $Q_{3,3}P_{1,4}-Q_{3,4}P_{1,3}=0$.
Finally, if also $Q_{3,4}\equiv0$ we get only a singularity with $z\neq 0$
at zeroes of $Q_{3,3}$
\qed

\roep Remark \num. The possibility $2\al$ cannot occur under the condition
$b_1+b_2=e_1+e_2+e_3+e_4-2$, because it implies $e_4>e_2-2$, contradicting
an inequality derived earlier.
\endroep
\endcomment

For $K3$ surfaces we need $b_1+b_2=e_1+e_2+e_3+e_4-2$.
We give a table listing the possibilities under this assumption,
cf.~\cite{\rfbra, Table A.1--A.4}. 

\comment
The solutions are found by deriving inequalities for the $e_i$.
E.g.~in case $\al$ we obtain from $b_1+b_2=e_1+e_2+e_3+e_4-2$ and
$b_1\leq e_1+e_3$ that $b_2\geq e_2+e_4-2$. Combined with $b_2\leq 2e_4$
this gives $e_2 \leq e_4+2$. From $b_1\leq 2e_2$ we get $e_1+e_3\leq
e_2+e_4+2$. Writing $e_1=e_2+k$, $e_2=e_4+l$ and $e_3=e_4+m$ we get 
$k+m\leq 2$ and $m \leq l \leq 2$.
\endcomment
The table lists the possible values for $(b_1,b_2)$ and gives
for each  pair the invariants $(e_1,e_2,e_3,e_4)$ of
the scrolls on which the curve can lie. These form one
deformation type  with adjacencies going vertically,
except $S_{e+2,e,e,e}$ and $S_{e+1,e+1,e+1,e-1}$ 
which do not deform into each other but are both deformations of 
$S_{e+2,e+1,e,e-1}$ and both deform to $S_{e+1,e+1,e,e}$.
Furthermore we give the number of moduli for each family.

In the table we also list the base locus of $|2H-b_1R|$
(which contains that of $|2H-b_2R|$).
The base locus is a subscroll, for which we use the following notation
\cite{\rfrch, 2.8}:
we denote by $B_a$
the subscroll corresponding to the subset of all $e_i$ with $e_i\leq a$,
defined by the equations $z^{\kl(j)}=0$ for $e_j>a$.
We give
the number and type of the singularities of the general
element; the number given in the second half
of \cite{\rfbra, Table A.2} is not correct. 

As example of the computations we look at
$(e+3,e,e-1,e-2)$ with $(b_1,b_2)=(2e,2e-2)$.
The two equations have the form
$$
\displaylines{
p_1xw+p_2xz+p_0y^2+p_3yx+p_6x^2 \cr
q_0z^2+q_0yw+q_3xw+q_1yz+q_4xz+q_2y^2+q_5yx+q_8x^2\;,
}
$$
where the index denotes the degree in $(s\cn t)$.
We first use coordinate transformations to simplify these
equations. By replacing $y$, $z$ and $w$ by suitable multiples
we may assume that the three constant polynomials are 1.
Now replacing $z$ by $z-{1\over2}q_1y-{1\over2}q_4x$ removes the $yz$ and $xz$
terms. 
We then replace $y$ by $y-q_3x$ to get rid of the $xw$ term. By changing
$w$ we finally achieve the form $z^2+yw+q_8x^2$.
By a change in $(s\cn t)$ we may assume that $p_1=s$. We now look at the
affine chart $(w=1,t=1)$ and find $y=-z^2-q_8x^2$, which we 
insert in the other equation to get an equation
of the form $x(s+p_2z+\cdots)+z^4$, which is an $A_3$.

We leave it again to the reader to analyse 
which further singularities can occur if $e_4=0$.

\section Infinitesimal deformations.

Deformations of cones over complete intersections on scrolls 
need not preserve the rolling factors format. We shall study in detail
those who do. Many deformations of negative weight are of this type.

\roep Definition \num.
A {\sl pure rolling factors  deformation\/} is 
a deformation in which the scroll is undeformed
and only the  equations on the scroll are perturbed. 
\endroep

This means that the deformation of the 
additional equations can be written with the rolling factors.
Such deformations are always unobstructed.
However this is not the only type of deformation for which the scroll
is not changed. In weight zero one can have deformations
inside the scroll, where the type $(b_1,\dots,b_l)$ changes.

\roep Definition \num.
A (general) {\sl rolling factors  deformation\/} is 
a deformation in which the scroll is deformed
and the  additional equations are written  in rolling factors
with respect to the deformed scroll.
\endroep

The equations for the total space of a 1-parameter 
rolling factors  deformation describe
a scroll of one dimension higher, containing
a subvariety of the same codimension,
again in rolling factors format. 
Deformations
over higher dimensional base spaces may be obstructed.
Again in weight zero one can have deformations
of the scroll, where also the type $(b_1,\dots,b_l)$ changes.

Finally there are {\sl non-scrollar deformations\/},
where the perturbation of the scrollar equations does not define
a deformation of the scroll.
Examples of this phenomenon are easy to find (but difficult
to describe explicitly). A trigonal canonical curve  is a divisor
in a scroll, whereas the general canonical curve of the same genus $g$
is not of this type: the codimension of the trigonal locus 
in moduli space is $g-4$.

\roep Example \num. To give an example of a deformation inside a scroll,
we let $C$ be a tetragonal curve in $\P^9$
with invariants $(2,2,2;3,1)$. Then there is a weight 0 deformation
to a curve of type  $(2,2,2;2,2)$. To be specific, let $C$ be given
by $P=sx^2+ty^2+(s+t)z^2$, $Q=t^3x^2+s^3y^2+(s^3-t^3)z^2$.
We do not deform the scroll, but only the additional equations:
$$
\vbox{\halign{
$#$\hfil&$#$\hfil\cr
x_0^2+y_0y_1+z_0^2&{}+\ep(z_1^2-x_1^2)  \cr
x_0x_1+y_1^2+z_0z_1&{}+\ep(z_1z_2-x_1x_2)  \cr
x_1^2+y_1y_2+z_1^2&{}+\ep(y_0y_1+z_0z_1)  \cr
x_1x_2+y_2^2+z_1z_2&{}+\ep(y_1^2+z_1^2)  \cr
\noalign{\medskip}
 x_1x_2+y_0^2+z_0^2-z_1z_2 \hidewidth\cr
 x_2^2+y_0y_1+z_0z_1-z_2^2\hidewidth\cr
}}
$$
For $\ep\neq0$ we can write the ideal as
$$
\vbox{\halign{
$#$\hfil&$#$\hfil\cr
x_0^2+y_0y_1+z_0^2&{}+\ep(z_1^2-x_1^2)  \cr
x_0x_1+y_1^2+z_0z_1&{}+\ep(z_1z_2-x_1x_2)  \cr
x_1^2+y_1y_2+z_1^2&{}+\ep(z_2^2-x_2^2)  \cr
\noalign{\medskip}
x_0x_1+y_1^2+z_0z_1&{}+\ep(y_0^2+z_0^2) \cr
x_1^2+y_1y_2+z_1^2&{}+\ep(y_0y_1+z_0z_1)  \cr
x_1x_2+y_2^2+z_1z_2&{}+\ep(y_1^2+z_1^2)  \cr
}}
$$
We can describe this deformation in the following way.
Write $Q=sQ_s+tQ_t$. The two times three equations above are
obtained by rolling factors from $sP-\ep Q_t$ and $tP+\ep Q_s$.
We may generalise this example.

\proclaim Lemma \num. Let $V$ be a complete intersection of divisors 
of type $aH-b_1R$, $aH-b_2R$, given by equations $P$, $Q$.
If it is possible to write $Q=sQ_s+t^{b_1-b_2-1}Q_t$ then
the equations $sP-\ep Q_t$, $t^{b_1-b_2-1}P+\ep Q_s$ give a deformation
to a complete intersection of type $aH-(b_1-1)R$, $aH-(b_2+1)R$.

In general one has to combine such a deformation with a 
deformation of the scroll.

\roep {\rm(\num)} Non-scrollar deformations.
\roep Example \num.
As mentioned before such deformations must exist in weight zero
for trigonal cones.
We proceed with the explicit computation of embedded deformations.
We start from the normal bundle exact sequence
$$
0 \lra  N_{S/C} \lra N_C \lra N_S\otimes \sier C \lra 0 \;.
$$
As $C$ is a curve of type $3H-(g-4)R$ on $S$ we have that
$C\cdot C=3g+6$ and $H^1(C,N_{S/C})=0$. So we are interested in 
$H^0(C, N_S\otimes \sier C)$, and more particularly in the
cokernel of the map $H^0(S,N_S)\lra H^0(C, N_S\otimes \sier C)$,
as $H^0(S,N_S)$ gives   deformations of the scroll.

\proclaim Proposition \num.
The
cokernel of the map $H^0(S,N_S)\lra H^0(C, N_S\otimes \sier C)$
has dimension $g-4$.

\roep Proof.
An element of $H^0(C, N_S\otimes \sier C)$ is a function 
$\ph$ on the equations
of the scroll such that the generators of the  module of 
relations map to zero in $\sier C$ and
it lies in the image of $H^0(S,N_S)$ if the function values can be lifted to
$\sier S$ such that the relations map to $0\in \sier S$.
Therefore we perform our computations in $\sier S$.

We have to introduce some more notation. Using
the equations described in (\labsect) 
\comment
In the notation introduced above, the scroll is
$$
\pmatrix{
x_0&x_1&\ldots&x_{e_1-1}&y_0&y_1&\ldots&y_{e_2-1}\cr
x_1&x_2&\ldots&x_{e_1  }&y_1&y_2&\ldots&y_{e_2}
}\;.
$$
The curve $C$ is given by a bihomogeneous equation
$$
F=A_{2e_1-e_2+2}x^3+B_{e_1+2}x^2y+C_{e_2+2}xy^2+D_{2e_2-e_1+2}y^3
$$
where $A_{2e_1-e_2+2}$ is a polynomial in $s$ and $t$ of degree  
$2e_1-e_2+2=3e_1-g+4$ and similarly for the other coefficients.
We have three types of equations,
\endcomment
we have three types of scrollar equations,
$f_{i,j}=x_ix_{j+1}-x_{i+1}x_{j}$, $g_{i,j}=y_iy_{j+1}-y_{i+1}y_{j}$
and mixed equations $h_{i,j}=x_iy_{j+1}-x_{i+1}y_{j}$.
The scrollar relations come from doubling a row in the matrix
and there are two ways to do this. The equations resulting from
doubling the top row can be divided by $s$, and the other ones by $t$,
so the result is the same.

A relation involving only equations of type $f_{i,j}$ gives the
condition 
$$
x s^{e_1-i-1}t^i\ph(f_{j,k})
-x s^{e_1-j-1}t^j\ph(f_{i,k})
+x s^{e_1-k-1}t^k\ph(f_{i,j}) =0 \in \sier C
$$
which may be divided by $x$. As the image $\ph(f_{i,j})$ is 
quadratic in $x$ and $y$ the resulting
left hand side cannot be a multiple of the equation of $C$, so we have
$$
 s^{e_1-i-1}t^i\ph(f_{j,k})
- s^{e-1-j-1}t^j\ph(f_{i,k})
+ s^{e_1-k-1}t^k\ph(f_{i,j}) =0 \in \sier S
$$
and the analogous equation involving only the $g_{i,j}$
equations.

For the mixed equations we get
$$
x s^{e_1-i-1}t^i\ph(h_{j,k})
-x s^{e_1-j-1}t^j\ph(h_{i,k})
+y s^{e_2-k-1}t^k\ph(f_{i,j}) = \psi_{i,j;k} P \in \sier S
$$
with $\psi_{i,j;k}$ of degree $e_1+e_2-3=g-5$
and  analogous ones involving $g_{i,j}$ with coefficients $\psi_{i;j,k}$.
These coefficients are not independent, but satisfy a systems of
equations coming from the syzygies between the relations. 
They can also be verified directly.
We obtain
$$
 s^{e_1-i-1}t^i\psi_{j,k;l}
- s^{e_1-j-1}t^j\psi_{i,k;l}
+ s^{e_1-k-1}t^k\psi_{i,j;l} =0 \in \sier S
$$
and 
$$
x s^{e_1-i-1}t^i\psi_{j;k,l}
-x s^{e_1-j-1}t^j\psi_{i;k,l}
+y s^{e_2-k-1}t^k\psi_{i,j;l} 
-y s^{e_2-l-1}t^l\psi_{i,j;k} 
=0 
$$
The last set of equations shows that $ s^{e_2-k-1}t^k\psi_{i,j;l} 
= s^{e_2-l-1}t^l\psi_{i,j;k} $ (rolling factors!) and therefore 
$\psi_{i,j;k}= s^{e_2-k-1}t^{k}\psi_{i,j;}$ with $\psi_{i,j;}$ of degree
$e_1-2$. This yields the equations
$$
 s^{e_1-i-1}t^i\psi_{j,k;}
- s^{e_1-j-1}t^j\psi_{i,k;}
+ s^{e_1-k-1}t^k\psi_{i,j;} =0
$$
Our next goal is to express all $\psi_{i,j;}$ in terms of 
the $\psi_{i,i+1;}$ (where $0\leq i \leq e_1-2$). 
First we observe by using the last equation for the
triples $(0,i,i+1)$ and $(i,i+1,e_1-1)$ that $\psi_{i,i+1;}$ is divisible
by $t^i$ and by $ s^{e_1-i-2}$ so $\psi_{i,i+1;}= s^{e_1-i-2}t^ic_i$ for some
constant $c_i$. By induction it then follows that
$\psi_{i,j;}= s^{e_1-i-2}t^ic_{j-1}+ s^{e_1-i-3}t^{i+1}c_{j-2}
+\cdots+ s^{e_1-j-1}t^{j-1}c_{i}$, so  the solution of the equations
depends on $e_1-1$ constants. Similarly one finds $e_2-1$ constants 
$d_i$ for the
$\psi_{i;j,k}$ so altogether $e_1+e_2-2=g-4$ constants. 

Finally we can solve for the perturbations of the equations. We give the
formulas in the case that all $d_i$ and all $c_i$ but one are zero, say 
$c_\ga=1$.
This implies that $\psi_{i,j;}=0$ if $\ga \notin [i,j)$ and
$\psi_{i,j;}= s^{e_1-i-j+\ga-1}t^{i+j-\ga-1}$ if $\ga \in [i,j)$; under the
last assumption $\psi_{i,j;k}= s^{g-4-i-j-k+\ga}t^{i+j+k-\ga-1}$.
We take $\ph(f_{i,j})=0$ if $\ga \notin [i,j)$. It follows that
for a fixed $k$ 
the $\ph(h_{i,k})$ with $i\leq \ga$ are related by rolling factors,
as are the $\ph(h_{i,k})$ with $i>\ga$. This reduces
the mixed equations with fixed $k$ to one, which can be solved
for in a uniform way for all $k$.
To this end we write the equation  $P$ as
$$
P=( s^{2e_1-e_2-\ga+2}A_\ga^{+}+t^{\ga+1}A^{-}_{2e_1-e_2-\ga+1})x^3
+y(B_{e_1+2}x^2+C_{e_2+2}xy+D_{2e_2-e_1+2}y^2)
$$
which we will abbreviate as $( s^{2e_1-e_2-\ga+2}A^{+}+t^{\ga+1}A^{-})x^3
+yE$.
We set 
$$

\matrix{
\ph(f_{i,j})&=&0,\hfill& $if$\ \ga \notin [i,j) \cr
\ph(f_{i,j})&=& s^{e_1-i-j-1+\ga}t^{i+j-\ga-1}E,\hfill& $if$\ \ga \in [i,j) \cr
\ph(g_{i,j})&=&0,\hfill&  \cr
\ph(h_{i,k})&=&- s^{e_2-1-k-i+\ga}t^{i+k}A^{-}x^2,
\hfill& $if$\ i\leq \ga \hfill\cr
\ph(h_{i,k})&=& s^{2e_1+1-i-k}t^{i+k-\ga-1}A^{+}x^2 ,
\hfill& $if$\  i> \ga \hfill\cr
}
$$
This is well defined, because all exponents of $ s$ and $t$ are positive.
\qed

A similar computation can be used to show that all elements of $T^1(\nu)$
with $\nu>0$ can be written rolling factors type. However,
even more is true, they can be represented as pure rolling factors
deformations, see \cite{\rfds}, where a direct argument is given.
\endroep

We generalise the above discussion to the case of a complete intersection
of divisors of type $aH-b_iR$ (with the same $a\geq2$) on a scroll
$$
\def\zz#1#2#3{z^{\kl(#1)}_{#2}&\ldots&z^{\kl(#1)}_{#3}}
\pmatrix{
\zz10{d_1-1}&\ldots&\zz k0{d_k-1}\cr
\zz11{d_1}&\ldots&\zz k1{d_k}\cr}.
$$
We have equations $f^{(\al\be)}_{ij}=z^{\kl(\al)}_iz^{\kl(\be)}_{j+1}
-z^{\kl(\al)}_{i+1}z^{\kl(\be)}_{j}$.
The lowest degree in which non rolling factors deformations can occur is
$a-3$. We get the conditions
$$
z^{\kl(\al)} s^{d_\al-i-1}t^i\ph\big(f^{(\be\ga)}_{jk}\big)
-z^{\kl(\be)} s^{d_\be-j-1}t^j\ph\big(f^{(\al\ga)}_{ik}\big)
+z^{\kl(\ga)} s^{d_\ga-k-1}t^k\ph\big(f^{(\al\be)}_{ij}\big) = 
\sum_l\psi^{(\al\be\ga)}_{ijk;l} P^\kl(l) 
$$
with the $\psi^{(\al\be\ga)}_{ijk;l}$ homogeneous 
polynomials in $( s\cn t)$ of degree $b_l-1$. 
The relations between these polynomials
come from the syzygies of the scroll: we add four of these
relations, multiplied with a term linear in the
$z^{\kl(\al)}$; then the left hand side becomes zero, leading to a
relation (in $\sier S$) between the $P^\kl(l)$. As we are dealing with 
a complete intersection, the relations are generated by Koszul relations.
Because the coefficients of the relation 
obtained are linear in the $z^{\kl(\al)}$,
they cannot lie in the ideal generated by the $P^\kl(l)$ 
(as $a\geq2$), so they
vanish and we obtain for each $l$ equations
$$
\displaylines{\qquad
z^{\kl(\al)} s^{d_\al-i-1}t^i\psi^{(\be\ga\de)}_{jkm;l} 
-z^{\kl(\be)} s^{d_\be-j-1}t^j\psi^{(\al\ga\de)}_{ikm;l}
\hfill\cr\hfill{} 
+z^{\kl(\ga)} s^{d_\ga-k-1}t^k\psi^{(\al\be\de)}_{ijm;l} 
-z^{\kl(\de)} s^{d_\de-m-1}t^m\psi^{(\al\be\ga)}_{ijk;l} 
=0\;.\qquad}
$$
Here some of the $\al$, \dots, $\de$ may coincide. If e.g. $\de$
is different from $\al$, $\be$ and $\ga$, then $\psi^{(\al\be\ga)}_{ijk;l} 
=0$. If there are at least four different indices (e.g. if the scroll
is nonsingular of dimension at least four) then $\de$ can always be chosen
in this way, so all coefficients vanish and every deformation
of degree $a-3$ is of rolling factors type.

Suppose now the scroll is a cone over a nonsingular 3-dimensional scroll,
i.e. we have three different indices at our disposal. Then every
$\psi^{(\al\be\ga)}_{ijk;l}$ with at most two different upper indices
vanishes, and the ones with three different indices satisfy rolling factors
equations.
We conclude that for pairwise different $\al$, $\be$, $\ga$
$$
\psi^{(\al\be\ga)}_{ijk;l}= s^{d-i-j-k-3}t^{i+j+k}\psi'_l
$$
with $d=d_\al+d_\be+d_\ga$ the degree of the scroll.

Finally, for the cone over a 2-dimensional scroll we get similar computations
as in the trigonal example above.

\proclaim Proposition \num.
\global\edef\labpropo{{\the\secno.\the\thmno}}
A tetragonal cone (with $g>5$) has non-scrollar deformations of degree $-1$
if and only if $b_2=0$. If the canonical curve lies on a Del Pezzo surface
then the dimension is $1$.
If the curve is bielliptic  then the dimension is $b_1=g-5$.

\roep Proof.
First suppose $e_3>0$. Then the only possibly non zero coefficients are the
$\psi'_l$, which have degree $b_l+2-\sum e_i$.
As $b_1+b_2=\sum e_i -2$ they do not vanish iff $b_2=0$.
In this case the computation yields one non rolling factors
deformation of the Del Pezzo surface on which the curve lies.

If $e_3=0$, then $b_2=0$. For a bielliptic curve the methods above
yield $(e_1-1)+(e_2-1)=b_1=g-5$ non-scrollar deformations
(a detailed computation is given in \cite{\rfbrap}).  
Suppose now that the curve lies on a (singular) Del Pezzo surface.
If $b_1=e_1>e_2=2$ then the equation $P$ contains the monomial $xz$
with nonzero coefficient, which we take to be $1$, while there is
no monomial $yz$. After a coordinate transformation we may assume 
that the same holds in case $e_1=e_2=2$. Let $\ph(h_{i,k})\equiv
\ze_{i,k}z \bmod (x,y)$. In the equation 
$$
x s^{e_1-i-1}t^i\ph(g_{j,k})
-y s^{e_2-j-1}t^j\ph(h_{i,k})
+y s^{e_2-k-1}t^k\ph(h_{i,j}) = \psi_{i;j,k} P 
$$
holding in $\sier S$
the monomial $yz$ occurs only on the left hand side, which
shows that the  $\ze_{i,k}$ are of rolling factors type in the
first index. Being constants,  they vanish.
This means that in the equation
$$
x s^{e_1-i-1}t^i\ph(h_{j,k})
-x s^{e_1-j-1}t^j\ph(h_{i,k})
+y s^{e_2-k-1}t^k\ph(f_{i,j}) = \psi_{i,j;k} P
$$
the monomial $xz$ does not occur on the left hand side
and therefore $\psi_{i,j;k}=0$. We find only $e_2-1=1$  
non rolling factors deformation. If $e_1=2$, $e_2=1$ we find one
deformation. Finally, if $e_1=3$, $e_2=1$ then there is only one
type of mixed equation. We have two constants $c_0$ and $c_1$.
Let the coefficient of $xz$ in $P$ be $p_0s+p_1t$.
We obtain the equations
$$
\displaylines{
s\ze_{1,0}-t\ze_{0,0}=c_0(p_0s+p_1t) \cr
s\ze_{2,0}-t\ze_{1,0}=c_1(p_0s+p_1t) \cr
}
$$
from which we conclude that $p_0c_0+p_1c_1=0$, giving again only one
non rolling factors deformation.
\qed

\comment
\roep Some remarks on $K3$-surfaces and Fano's.

Remark:
Relevant: results of Donagi-Morrison on hyperplane sections of
$K3$ surfaces.
Relevant: CLM on Fano's.
\endcomment

\roep {\rm(\num)} Rolling factors deformations of degree $-1$.

\noindent
We look at the miniversal deformation of the scroll:
$$
\pmatrix{
z^{\kl(1)}_0&\ldots&z^{\kl(1)}_{d_1-2}&z^{\kl(1)}_{d_1-1}&z^{\kl(2)}_0 &
 \ldots&z^{\kl(k)}_{d_k-2}&z^{\kl(k)}_{d_k-1}\cr
z^{\kl(1)}_1+\ze^{\kl(1)}_1&\ldots&z^{\kl(1)}_{d_1-1}+\ze^{\kl(1)}_{d_1-1}
&z^{\kl(1)}_{d_1}&z^{\kl(2)}_1 &
 \ldots&z^{\kl(k)}_{d_k-1}+\ze^{\kl(k)}_{d_k-1}&z^{\kl(k)}_{d_k}}
$$
To compute which of those deformations can be lifted to deformations
of a complete intersection on the scroll we have to compute perturbations
of the additional equations.

We assume that we have  a complete intersection
of divisors of type $aH-b_iR$ (with the same $a\geq2$).

Extending the notation introduced before we write the columns
in the matrix  symbolically as
$(z_\al,z_{\al+1}+\ze_{\al+1})$.
In order that this makes sense for all columns  we introduce dummy variables
$\ze^{\kl(i)}_0$ and $\ze^{\kl(i)}_{d_i}$ with the value $0$.

The Koszul type relations give no new conditions, but 
the relation
$$
P_{m+1}z_\be-P_mz_{\be+1}-\sum_\al p_{\al,m}f_{\be\al}=0
$$
gives as equation in the local ring for the perturbations 
$P_m'$ of $P_m$:
$$
P_{m+1}'z_\be-P_m'z_{\be+1}-\sum_\al p_{\al,m}
(\ze_{\al+1}z_\be-z_\al\ze_{\be+1})=0\;.
$$
In particular we see that we can look at one equation on the
scroll at a time. 
As $\sum p_{\al,m}z_\al=P_m$ the coefficient of $\ze_{\be+1}$
vanishes.
Because  $tz_\be-sz_{\be+1}=0$ we get a condition which is independent
of $\be$:
$$
sP_{m+1}'-tP_m'-s\sum_\al p_{\al,m} \ze_{\al+1}=0
$$
This has to hold in the local ring, but as the degree of the $ p_{\al,m}$
is lower than that of the equations defining the complete
intersection on the scroll (here we use the assumption that all degrees $a$
are equal), it holds on the scroll.
From it we derive the equation
\def\labeq{(S)}
$$
s^bP_{b}'-t^bP_0'=
 \sum_{m=0}^{b-1}\sum_\al s^{m+1}t^{b-m-1}p_{\al,m} \ze_{\al+1} \eqno \labeq
$$
which has to be solved with $P_b'$ and $P_0'$ polynomials in the
$z_\al$ of degree $a-1$. We determine the monomials on the right hand side.

The result depends on the chosen equations, but only on $P_0$ and $P_b$ 
and not on the intermediate ones, provided they are obtained by rolling
factors. 
\roep Example \num. 
\edef\labex{{\the\secno.\the\thmno}}
Let $b=4$. We take variables $y_i=s^{3-i}t^iy$,
$z_i=s^{3-i}t^iz$ with deformations $\eta_i$, $\ze_i$, and roll from 
$y_0z_0$ to $y_2z_2$ in two different ways:
$$
\displaylines{
y_0z_0 \to y_1z_0 \to y_1z_1 \to y_2z_1 \to y_2z_2 \cr
y_0z_0 \to y_0z_1 \to y_0z_2 \to y_1z_2 \to y_2z_2 }
$$
This gives  as right-hand side of the equation $\labeq$ in the two cases
$$
\displaylines{
s^4t^3z\eta_1 + s^4t^3 y\ze_1 +s^5t^2 z\eta_2 +s^5t^2 y\ze_2 \cr
s^4t^3y\ze_1 + s^5t^2 y\ze_2 +s^4t^3 z\eta_1 +s^5t^2 z\eta_2  }
$$
which is the same expression.
Similarly, if we roll from $z_0^2$ to $z_2^2$ we get 
$$
2s^4t^3z\ze_1 + 2s^5t^2 z\ze_2 
$$
However, if we roll in the last step from $y_1z_2$ to $y_1z_3$
we get 
$$
s^4t^3y\ze_1 + s^5t^2 y\ze_2 +s^4t^3 z\eta_1
$$
(remember that we have no deformation parameter $\ze_3$).
\endroep

To analyse the general situation it is convenient to use 
multi-index notation. 
The equation $P$ of a divisor in $|aH-bR|$ may then be written as
$$
P = \sum_{|I|=a}\sum_{j=0}^{\ip<e,I>-b} p_{I,j}s^{\ip<e,I>-b-j}t^jz^I\;.
$$
Here $e=(e_1,\dots,e_k)$ is the vector of
degrees and $z^I$ stands for 
$(z^{\kl(1)})^{i_1}\cdot \dots \cdot(z^{\kl(k)})^{i_k}$.

\proclaim Proposition \num.
The lifting condition for the equations $P_m$ is that 
for each $I$ with $|I|=a-1$ and $\ip<e,I><b-1$ the following 
$b-\ip<e,I>-1$
linear equations hold:
$$
\sum_{l=1}^k\sum_{j=0}^{\ip<e,I+\de_l>-b}(i_l+1)
p_{I+\de_l,j}\ze^{\kl(l)}_{j+n}=0\;,
$$
where $0<n<b-\ip<e,I>$

\roep Proof.
We look at a monomial $s^{\ip<e,I'>-b-j}t^jz^{I'}$. In rolling from
$P_0$ to $P_m$ we go from $z_A$ to $z_{A+B}$. Here we write
a monomial as product of $a$ factors: $z_{\al_1}\cdots z_{\al_a}$
with $i'_l$ factors of type $l$.
Let $I'=I+\de_l$ with $\de_l$ the $l$th unit vector.
The monomial leads to an expression in which the coefficient of
$z^{I}$ is 
$$
\sum_{\{q\mid\al_q\ {\rm of\ type}\ l\}}\sum_{r=1}^{\be_q}  
s^{\ip<e,I>+r-j+\al_q}t^{b-r+j-\al_q}\ze_{\al_q+r}
$$
We stress that the choice of $\al_q$ can be very different for
different $j$.

We collect all contributions and look at the coefficient
of $s^{\ip<e,I>+n}t^{b-n}z^{I}$ with $0<n<b-\ip<e,I>$. This cannot
be realised as left-hand side of equation $\labeq$.
Because $b-n=b-r+j-\al_q$ this coefficient is
$$
\sum_{l=1}^k\sum_{j=0}^{\ip<e,I+\de_l>-b}(i_l+1)
p_{I+\de_l,j}\ze^{\kl(l)}_{j+n}=0\;,
$$
We note that all terms really occur:  in rolling from
$z_A$ to $z_{A+B}$ we have to increase the  $q$th factor
sufficiently many times, because $\ip<e,I><b-1$. 
\qed

\roep Example \num: trigonal cones.
Let the curve be given by the bihomogeneous equation
$$
F=A_{2a-m+2}z^3+B_{a+2}z^2w+C_{m+2}zw^2+D_{2m-a+2}w^3
$$
then there are only conditions for $I=(0,2)$, i.e.~for $w^2$,
as $a+m=g-2>b-1$. So if $2m<b-1=g-5$ we get $b-2m-1=a-m-3$
equations on the deformation variables $\ze_1$, \dots, $\ze_{a-1}$,
$\om_1$, \dots, $\om_{m-1}$
$$
\sum_{j=0}^{m+2}c_{j}\ze_{j+n}
+3\sum_{j=0}^{2m-a+2}d_{j}\om_{j+n}=0\;,
$$
as stated in \cite{\rfds, 3.11}.
We have a system of linear equations so we can write the coefficient
matrix. It consists of two blocks $(\, {\cal C}\mid{\cal D}\,)$ with
${\cal C}$ of the form
$$
\pmatrix{
c_0 & c_1 & c_2 & \dots & c_{m+2}& 0 & 0& \dots &0&0&0  \cr
0   & c_0 & c_1 & \dots & c_{m+1}& c_{m+2}& 0&\dots &0&0&0 \cr
0   & 0   & c_0 & \dots & c_m & c_{m+1}& c_{m+2}&\dots &0&0&0  \cr
    &     &     & \ddots&     &        &        & \ddots\cr
0   & 0   & 0 & \dots & c_0 & c_1 & c_2 & \dots & c_{m+2}&0&0     \cr
0   & 0   & 0 & \dots & 0  & c_0 & c_1 & \dots & c_{m+1}& c_{m+2}&0 \cr
0   & 0   & 0 & \dots & 0& 0& c_0 &  \dots & c_m & c_{m+1}& c_{m+2} \cr
}
$$
and $\cal D$ similarly.
Obviously this system has maximal rank.
\endroep

The Proposition gives  a system of linear equations
and we call the coefficient matrix {\sl lifting matrix\/}.
It was introduced for tetragonal cones in \cite{\rfbrap}.

In general the lifting matrix will have maximal rank, 
but it is a difficult question to decide when this happens.

\roep Example \num: trigonal $K3$s.
We take the invariants $(e,e,e)$, $b=3e-2$ with $e\geq 3$. The $K3$ lies on 
$\P^2\times\P^1$ and is given by an equation  
of bidegree $(3,2)$. Now there are six $I$s with $|I|=2$ each giving rise
to $e-3$ equations in $3(e-1)$ deformation variables.
In general the matrix has maximal rank, but for special surfaces
the rank can drop.
Consider an equation of type 
$p_1x^3+p_2y^3+p_3z^3$ with the $p_i$ quadratic polynomials in $(s\cn t)$
without common or multiple zeroes. Then the surface is smooth.
The lifting equations corresponding to the quadratic monomials
$xy$, $xz$ and $yz$ vanish identically and the lifting matrix reduces to
a block-diagonal matrix of rank $3(e-3)$. The kernel has dimension $6$,
but the corresponding deformations are  obstructed: an
extension of the $K3$ would be a Fano 3-fold
with isolated singularities lying as divisor of type $3H-bR$ on a scroll
$S(e_1,e_2,e_3,e_4)$, and a computation  reveals
that such a Fano can only exist for $\sum e_i \leq 8$.
\endroep

We can say something more for the lifting conditions
coming from one quadratic equation.

\proclaim Proposition \num. 
The lifting matrix for one quadratic equation has dependent rows 
if and only if the 
generic fibre has a singular point on the subscroll $B_{b-1}$.

\roep Proof.
The equation $P$ on the scroll can be written in the form 
${}^tz\Pi z$ with $\Pi$ a symmetric $k\times k$
matrix with  polynomials
in $(s\cn t)$ as entries. The condition that there is a singular section 
of the form  $z=(0,\dots,0,z^{\kl(l+1)}(s,t),\dots,z^{\kl(k)}(s,t))$
with $e_l>b-1$
is that $ {}^tz\Pi=0$  or $ {}^tz_{>l}\Pi_{>l}=0$
where $z_{>l}=(z^{\kl(l+1)}(s,t),\dots,z^{\kl(k)}(s,t))$ and $\Pi_{>l}$
is the matrix consisting of the last $k-l$ rows of $\Pi$.
The resulting system of equations for
the coefficients of the polynomials $z^{\kl(i)}(s,t)$
gives exactly the lifting matrix.
\qed

\roep {\rm(\num)} Tetragonal curves.
Most of the following results are contained in the preprint \cite{\rfbrap}.
We have two equations on the scroll and the lifting
matrix $M$ can have rows coming from both equations.
We first suppose that $b_2>0$. Then also $e_3>0$ and the number of
columns of $M$ is always $\sum (e_i-1)=g-6$, but the number of rows depends
on the values of  $(e_1,e_2,e_3;b_1,b_2)$: 
it is $\sum_{i,j}\max(0,b_i-e_j-1)$.

\proclaim Theorem \num.
Let $X$ be the cone over a tetragonal canonical curve and suppose that $b_2>0$.
Then $\dim T^1_X(-2)=0$.
Suppose that  the $g^1_4$ is not composed with an involution
of genus $\frac{b_2}2+1$.
\item{1)} If $b_1<e_1+1$ or $b_2<e_3+1$ or $g\leq 15$
then $\dim T^1_X(-1)=9 +\dim \Cork M$. 
\item{2)} If $b_1\geq e_1+1$, $b_2\geq e_3+1$ and  $g> 15$ then
$ 9 +\dim \Cork M \leq \dim T^1_X(-1) \leq \frac{g+3}6+6 +\dim \Cork M$
and the maximum is obtained for
$g$ of the form $6n-3$ and $(e_1,e_2,e_3;b_1,b_2)=(3n-2,2n-2,n-2;4n-4,2n-4)$.
\item{3)} For generic values of the moduli $\dim \Cork M=0$.

\roep Proof. 
If $b_2>0$ there are only rolling factors deformations in 
negative degrees. In particular $\dim T^1_X(-2)=0$.
The number of pure rolling factors deformations
is $\rho=\sum_{i,j}\max(e_j-b_i+1,0)$.
The number of rows in the lifting matrix is $\sum_{i,j}\max(0,b_i-e_j-1)=
3(b_1+b_2)-2(e_1+e_2+e_3+3)+\rho=g-15+\rho$. If $\rho>9$ the number of
rows exceeds the number of columns and $\dim T^1_X(-1)=\rho +\dim \Cork M$;
otherwise it is $9+\dim \Cork M$. So we have to estimate $\rho$.

As the $g^1_4$ is not composed we have  
$b_1\leq e_1+e_3$. Together with  $b_1\leq 2e_2$ we get  $3b_1\leq 2g-6$
and $3b_2\geq g-9$; from $e_1\leq \frac{g-1}2$ we now
derive $e_1-b_2+1\leq {g-1\over2}+1-{g-9\over 3}$.
Also 
$b_2=e_1+e_2+e_3-2-b_1\geq e_2-2$, so $e_2-b_2+1\leq 3$.

2) Suppose first that $b_1\geq e_1+1$ and $b_2\geq e_3+1$.
Then $\rho = \max(0,e_1-b_2+1)+\max(0,e_2-b_2+1)\leq {g+3\over 6}+6$.
Equality is achieved iff $e_1=(g-1)/2$, $b_2=(g-9)/3$ and
$e_2=b_2+2$, so $g$ 
has the form $6n-3$ and $(e_1,e_2,e_3;b_1,b_2)=(3n-2,2n-2,n-2;4n-4,2n-4)$.

1) In all other cases $\rho\leq 9$: 
if  $b_1\geq e_1+1$, but $b_2< e_3+1$ then $\rho= 
(e_1-b_2+1)+(e_2-b_2+1)+(e_3-b_2+1)=g-3b_2\leq 9$.
If $e_2+1\leq b_1 <e_1+1$ then $\rho =
(e_1-b_1+1)+(e_1-b_2+1)+ \max(0,e_2-b_2+1)+\max(0,e_3-b_2+1)=
2e_1+7-g+\max(0,e_2-b_2+1)+\max(0,e_3-b_2+1)$.
As $b_2>0$ we have that $\max(0,e_2-b_2+1)+\max(0,e_3-b_2+1)
>\max(0,e_2+e_3-b_2+1)$. But from $b_1\leq e_1$ it follows that
$b_2\geq e_2+e_3-2$. So $\rho \leq (g-1)+7-g+3=9$.
If $b_1<e_2+1$ then $\rho\leq 2e_1+2e_2+4-2b_1-2b_2+2e_3=8$.

3) It is easy to construct lifting matrices of maximal rank
for all possible numbers of blocks occurring. 
\qed

Now we consider the case that the $g^1_4$ is composed with an involution
of genus $g'=\frac{b_2}2+1$. So if $b_2>0$, then $g'>1$.
After a coordinate transformation
we may assume that the surface $Y$ is singular along the section $x=y=0$,
so its equation depends only on $x$ and $y$: $P=P(x,y;s,t)$.
We may assume that $Q$ has the form $Q=z^2+Q'(x,y;s,t)$.
Let $M_{xy}$ be the submatrix of the lifting matrix consisting of the
blocks coming from $P$ and $Q'$ and the $\xi$ and $\eta$ deformations.

\proclaim Theorem \num.
Let $X$ be  a tetragonal canonical cone 
such that the $g^1_4$ is composed with an involution
of genus $g'>1$.
Then $\dim T^1_X(-1) = e_1+e_2-2e_3+6+\Cork M_{xy}$.

\roep Proof.
The rows in the lifting matrix $M$ coming from the first equation
and the variable $z$ vanish identically. The second equation
gives a $z$-block which is an identity matrix
of size $b_2-e_3-1=e_3-1$, so all $\zeta$
variables have to vanish. What remains is the matrix
$M_{xy}$ which has $e_1+e_2-2$ columns.
The number of rows is $\max(0,e_2-e_3-3)+\max(0,e_1-e_3-3)+
\max(0,2e_3-e_1-1)+\max(0,2e_3-e_2-1)$. We estimate the last two terms
with $e_3-1$ and the  first two by $e_2-e_3$, resp.~$e_1-e_3$.
Therefore the number of rows is at most $e_1+e_2-2$. For each term
which contributes $0$ to the sum we have pure rolling
factors deformations, so if the matrix has maximal
rank the dimension of $T^1_X(-1)$ is $e_1+e_2-2-(2e_3-8)$.
\qed

\roep Example \num.
It is possible that the lifting matrix $M$ does not have full
rank even if the $g^1_4$ is not composed.
An example with  invariants $(6,5,5;7,7)$ is the curve given by
the equations
$(s^5+t^5)x^2+s^3y^2+t^3z^2$, $s^5x^2+(s^3-t^3)(y-z)^2+2t^3z^2$.
The matrix is 
$$
\pmatrix{
0&\dots &0&|& 2 & 0 & 0 & 0&|& 0 & 0 & 0 & 0 \cr
0&\dots &0&|& 0 & 0 & 0 & 0&|& 0 & 0 & 0 & 2 \cr
\noalign{\smallskip\hrule}
\noalign{\smallskip}
0&\dots &0&|& 2 & 0 & 0 & -2&|& -2 & 0 & 0 & 2 \cr
0&\dots &0&|& -2 & 0 & 0 & 2&|& 2 & 0 & 0 & 2 \cr
}\;.
$$
\endroep

Finally we mention the case $b_2=0$. Bielliptic curves
($e_3=0$) are treated in \cite{\rfcm},  curves on a Del Pezzo
in \cite{\rfbrap} (but he overlooks those with $e_3=0$). Now there
is only one equation coming from $Q$, which can be perturbed arbitrarily.
As the $z$ variable does not enter the scroll, we have one coordinate
transformation left. The lifting matrix involves only rows coming
from the equation $P$. One checks that
the matrix $M$ resp.~$M_{xy}$ has maximal rank
and the number of rows does not exceed the number of columns.
Together with the number of non-scrollar deformations
(Prop.~\labpropo) this yields the following result, where we have
excluded the complete intersection case $g=5$.

\proclaim Proposition \num.
Let $X$ be the cone over a tetragonal canonical curve $C$ with $b_2=0$
and $g>5$.
Then $\dim T^1_X(-2)=1$.
\item{1)} If $C$ lies on a Del Pezzo surface
then $\dim T^1_X(-1)=10$. 
\item{2)} If $C$ is bielliptic ($e_3=0$), then $\dim T^1_X(-1)=2g-2$.

\roep Remark \num.
For all non-hyperelliptic canonical cones the dimension
of $T^1_X(\nu)$ with $\nu\geq0$ is the same. The
Wahl map easily gives 
$\dim T^1_X(0)= 3g-3$, $\dim T^1_X(1)= g$, $\dim T^1_X(2)= 1$
and $\dim T^1_X(\nu)=0$ for $\nu\geq3$
(see e.g.~\cite{\rfds}, 3.3).

\section Rolling factors obstructions.

Rolling factors deformations can be obstructed. We first give a general
result on the dimension of $T^2$. For the case of quadratic equations 
on the scroll one can actually write down the base equations.

\proclaim Proposition \num.
Let $X$ be the cone over a complete intersection
of divisors of type $aH-b_iR$ with $b_i>0$ (and the same $a\geq2$)
on a scroll.
If $a>2$, then $\dim \ttw(-a)=\sum (b_i-1)$, and 
$\dim \ttw(-a)\geq \sum (b_i-1)$ in
case $a=2$.

\roep Proof.
Let $\psi\in \Hom(R/R_0,\sier X)$ be an homogeneous element of
degree $-a$. 
The degree of $\psi(R_{\al,\be,\ga})$ is $3-a$, so $\psi$ vanishes on
the scrollar relations, if $a>2$. If $a=2$ we can assert that the
functions vanishing on the scrollar relations span a subspace of 
$\ttw(-2)$.

As the degree of the relation $R^n_{\al,m}$ is $a+1$,
the image $\psi(R^n_{\al,m})$ is a linear function of the coordinates.
The relations
$$
R^n_{\al,m}z_\be-R^n_{\be,m}z_{\al}-
\sum R^n_{j,k,\ga}p^n_{\ga,m}=
P^\kl(n)_m f_{\al,\be}-f_{\al,\be}P^\kl(n)_m.
$$
imply that the $\psi(R^n_{\al,m})$ are also in rolling factors form.
A basis (of the relevant subspace) of
$\Hom(R/R_0,\sier X)(-a)$ consists of the $2\sum b_i$ elements
$\psi^i_{l,s}(R^j_{\al,m})=\delta_{ij}\delta_{lm}z_\al$,
$\psi^i_{l,t}(R^j_{\al,m})=\delta_{ij}\delta_{lm}z_{\al+1}$, 
where $0\leq m <b_j$. The
image of $P^\kl(i)_m$ in $\Hom(R/R_0,\sier X)(-a)$ is $\psi^i_{m-1,s} -
\psi^i_{m,t}$, if $0<m<b_i$, $-\psi^i_{0,t}$ for $m=0$, and
$\psi^i_{b_i-1,s}$ for $m=b_i$. The quotient has dimension $\sum (b_i-1)$. 
\qed

For $a=2$  only the rolling factors obstructions 
will contribute to the base equations.
A more detailed study could reveal if there
are other obstructions. Typically this can happen, if there exist 
non-scrollar deformations.
As example we mention Wahl's result for tetragonal cones 
that $\dim\ttw(-2) =g-7=b_1+b_2-2$, if $b_2>0$, whereas 
for a curve on a Del Pezzo the dimension is $2(g-6)$
\cite{\rfwasq, Thm.~5.9}.

In the quadratic case we can easily write the base equations,
given a first order lift of the scrollar deformations.
We can consider
each equation on the scroll separately, so we will suppress  
the upper index of the additional equations in our notation.
\comment
Recall the assumption that 
two consecutive equations are of the form
$$
\eqalign{
P_m=&\sum_\al p_{\al,m}z_\al,\cr
P_{m+1}=&\sum_\al p_{\al,m}z_{\al+1}.\cr}
$$
\endcomment
We may assume that we have pure rolling factors deformations
$\rho_\al$ and that the lifting conditions are satisfied.
We can write the perturbation of the equation $P_m$
as 
$$
P_m(z)+P_m'(z,\ze,\rho)\;.
$$
Note that $P'_m$ is linear in $z$. Now we have the following
result \cite{\rfsth}.

\proclaim Proposition \num.
The maximal extension of the
infinitesimal deformation defined by the $P'_m$
is given by the $b-1$ base equations
$$
P_m'(\ze,\ze,\rho)-P_m(\ze)=0\;,
$$ 
with $1\leq m\leq b-1$.

\roep Proof.
We also suppress $\rho$ from the notation.
We have to lift the relations $R_{\be,m}$. 
As the lifting equations are satisfied we can write
$$
P_{m+1}'(z,\ze)z_\beta-P_m'(z,\ze)z_{\beta+1}
-\sum p_{\alpha,m}(z)z_{\alpha+1}z_\be =
\sum f_{\beta\gamma}d_{\ga}(\ze)\;,
$$
because the left hand side lies in the ideal of the scroll.
This identity involving quadratic monomials in the $z$-variables
can be lifted to the deformation of the scroll.
We write $\tf_{\beta\alpha}$
for the deformed equation $(z_\beta+\ze_\beta)z_{\alpha+1}-
z_{\beta+1}(z_\alpha+\ze_\alpha)$.
We get
$$
P_{m+1}'(z+\ze,\ze)z_\beta-P_m'(z,\ze)(z_{\beta+1}+\ze_{\be+1})
-\sum p_{\alpha,m}(z+\ze)z_{\alpha+1}z_\be =
\sum \tf_{\beta\gamma}d_\ga(\ze)\;.
$$
We now lift the relation  $R_{\beta,m}$: 
$$
\displaylines{\quad
\left(
   P_{m+1}(z) +P_{m+1}'(z+\ze,\ze)-P_{m+1}(s)
\right)
  z_\beta-
\left(
  P_m(z)+P_m'(z,\ze)
\right)
(z_{\beta+1}+\ze_{\be+1})
\hfill\cr
\hfill{} -\sum \tf_{\beta\alpha}p_{\alpha,m}(z)
-\sum \tf_{\beta\gamma}d_\ga(\ze)=
0.\quad\cr}
$$ 
If $1\leq m\leq b-1$, then
$P_m$ occurs in a relation as first and as second term.
Therefore  $P_m'(z,\ze)$ 
and $P_m'(z+\ze,\ze)-P_m(\ze)$ have to be equal. These
equations correspond to the $b-1$ elements of $\ttw(-2)$,
constructed above.   \qed

\roep Example \num. 
We continue with our rolling factors example \labex.
We look at two ways of rolling:
$$
\displaylines{
y_0z_0 \to y_1z_0 \to y_1z_1 \to y_2z_1 \to y_2z_2 \cr
y_0z_0 \to y_0z_1 \to y_0z_2 \to y_1z_2 \to y_1z_3 }
$$
The equation for $P_0'$ and $P_4'$ has a unique solution
with $P_0'=0$.
We get
$$
\matrix{
P_0' &=&\hfill  0,                     &   \hfill        0 \cr
P_1' &=&\hfill\eta_1 z_0 ,             &  \hfill    y_0 \ze_1  \cr
P_2' &=&\hfill \eta_1 z_1+y_1\ze_1  ,  &  \hfill y_1\ze_1 + y_0\ze_2 \cr
P_3' &=&\hfill \eta_1 z_2+y_2\ze_1 +\eta_2z_1 ,
                    & \hfill y_2\ze_1+y_1\ze_2+\eta_1z_2 \cr
P_4' &=&\hfill\eta_1z_3+y_3\ze_1+\eta_2z_2+y_2\ze_2,
                                 &\hfill y_3\ze_1+y_2\ze_2+\eta_1z_3}
$$
The resulting base equations are in both cases
$$
0, \eta_1\ze_1, \eta_1\ze_2+\eta_2\ze_1
$$
\endroep

In general the quadratic base equations are not uniquely determined. 
They can be modified by
multiples  of the linear lifting equations, if such are present. 
The other source of
non-uniqueness is the possibility of coordinate transformations using
the pure rolling factors variables. 

\proclaim Theorem \num.
Let $P=\sum p_{I,k}s^{\ip<e,I>-b-k}t^kz^I$ define a divisor
of type $2H-bR$. It leads to quadratic base equations $\pi_1$, \dots,
$\pi_{b-1}$. The  coefficient  $p_{I,k}$ gives the following  term in $\pi_m$.
We write $z^I=xy$ and assume that $e_x\geq e_y$.
\item{I.} If $e_x<b$ then for $m\leq k$ the term  is
         $-\sum_{l=m}^k \eta_{k-l+m}\xi_{l}$, while for $m>k$ it is
         $$
         \sum_{l=\max (k+m-e_y+1,k+1)}^{\min(e_x-1,m-1)}
         \hskip-2em\eta_{k-l+m}\xi_{l}\;.
         $$
\item{II.} If $e_x\geq b$ then for $m\leq k+b-e_x$ the term is
         $-\sum_{l=m+e_x-b}^k \eta_{k-l+m}\xi_{l}$, while for 
         $m> k+b-e_x$ it is
         $$
         \sum_{l=\max (k+m-e_y+1,k+1)}^{\min(e_x-b+m-1,k+m-1)}
         \hskip-2em\eta_{k-l+m}\xi_{l}\;.
         $$
Furthermore, if $e_x\geq b$ the $e_x-b+1$  pure rolling factors
deformations involving $x$ contribute 
$\rho_0 \xi_m+\cdots +\rho_{e_x-b}\xi_{m+e_x-b}$
to  $\pi_m$.

\edef\labthm{{\the\secno.\the\thmno}}
\roep Proof.
We have to choose explicit equations $P_m$. The monomial 
$s^{e-x+e_y-b-k}t^kxy$ gives a rolling monomial $x_{i(m)}y_{j(m)}$, where
$i(m)+j(m)=k+m$. Let $i(0)=i$, $i(b)=i'$, $j(0)=j$ and $j(b)=j'$.
We have to compute $P_m'$.
Equation $\labeq$ gives  
$$
s^bP_{b}'-t^bP_0'=
 \sum_{l=1}^{j'-j} s^{e_x-k+j+l}t^{k+b-j-l}x\eta_{j+l}
+ \sum_{n=1}^{i'-i} s^{e_y-k+i+n}t^{k+b-i-n}y\xi_{i+n}\;,
$$
which we rewrite as
$$
s^bP_{b}'-t^bP_0'=
 \sum_{l=j+1}^{j'} s^{e_x-k+l}t^{k+b-l}x\eta_{l}
+ \sum_{l=i+1}^{i'} s^{e_y-k+l}t^{k+b-l}y\xi_{l} \;.
$$

\roep Case I: $e_x< b$.
The condition $k+b\leq e_x+e_y$ implies $k< e_y$.
We solve for $P_0'$:
$$
P_0'=
    -\sum_{l=j+1}^{k} x_{k-l}\eta_{l}
    -\sum_{l=i+1}^{k} y_{k-l}\xi_{l}   \;.
$$
For the $P_m'$ we formally write the formula
$$
P_m'=
 -\sum_{l=j+1}^{k} x_{k-l+m}\eta_{l}
- \sum_{l=i+1}^{k} y_{k-l+m}\xi_{l}
+\sum_{l=j+1}^{j(m)} x_{k-l+m}\eta_{l}
+ \sum_{l=i+1}^{i(m)} y_{k-l+m}\xi_{l} \;.
$$
This expression can involve non-existing $x$ or $y$ variables:
for $y$ this happens if $k-l+m>e_y$, or $l<m+k-e_y$. The terms in the
two sums involving $y$ cancel. If $i(m)<k$, then the smallest 
non-cancelling term has $l=i(m)+1$ and $i(m)+1\geq i(m)+j(m)-e_y
=k+m-e_y$. If $i(m)>k$ we a sum of positive terms starting with
$k+1$. If $k<l<m+k-e_y$ then our monomial contributes to the
lifting conditions, and we can leave out this term. The sum
therefore now starts at $\max(k+1,m+k-e_y)$. Keeping this in mind
we determine the term in the base equation $\pi_m$ from the formal formula.
To this end we change the summation variable in the sums 
containing $x$-variables and arrive, using $i(m)+j(m)=k+m$,  at
$$
\displaylines{\quad
 -\sum_{l=m}^{m+i-1} \xi_{l}\eta_{k-l+m}
- \sum_{l=i+1}^{k} \eta_{k-l+m}\xi_{l}
+\sum_{l=i(m)}^{m+i-1} \xi_{l}\eta_{k-l+m}
+ \sum_{l=i+1}^{i(m)} \eta_{k-l+m}\xi_{l}-\xi_{i(m)}\eta_{j(m)} 
\hfill\cr\hfill{}=
 -\sum_{l=m}^{m+i-1} \xi_{l}\eta_{k-l+m}
- \sum_{l=i+1}^{k} \eta_{k-l+m}\xi_{l}
+ \sum_{l=i+1}^{m+i-1} \eta_{k-l+m}\xi_{l}\;.
\qquad
}
$$
If $k\geq m$ the terms from $l=m$ to $l=k$ occur twice with
a minus sign and once with a plus. Otherwise all negative terms cancel,
but we have to take the lifting conditions into account. 

\roep Case II: $e_x\geq b$.
Now there are $e_x-b+1$ pure rolling factors deformations present:
we can perturb $P_m$ with $\rho_0 x_m+\cdots +\rho_{e_x-b}x_{m+e_x-b}$.
These contribute $\rho_0 \xi_m+\cdots +\rho_{e_x-b}\xi_{m+e_x-b}$
to the equation $\pi_m$.

We can roll using only the $x$ variable: 
$x_{i+m}y_{j}$, with $i+j=k$ and $i+b\leq e_x$. We
take $i=k$ if $k+b\leq e_x$ and $i=e_x-b$ otherwise.
We get
$$
s^bP_{b}'-t^bP_0'=
\sum_{l=i+1}^{i+b} s^{e_y-k+l}t^{k+b-l}y\xi_{l} \;.
$$
We solve:
$$
P_0'=
    -\sum_{l=i+1}^{k} py_{k-l}\xi_{l}   
$$
and
$$
P_m'=
- \sum_{l=i+1}^{k} y_{k-l+m}\xi_{l}
+ \sum_{l=i+1}^{i+m} y_{k-l+m}\xi_{l} \;.
$$
Again if $k<l<m+k-e_y$  our monomial contributes to the
lifting conditions, and the sum
starts at $\max(k+1,m+k-e_y)$. We get as contribution to $\pi_m$
$$
- \sum_{l=i+1}^{k} \eta_{k-l+m}\xi_{l}
+ \sum_{l=i+1}^{i+m-1} \eta_{k-l+m}\xi_{l}
\;.
$$
Taking the lifting conditions and our choice of $i$
into account we get the statement
of the theorem.
\qed

\roep Example \num: Case I. 
Let $b=7$, $e_x=5$ and $e_y=4$. Consider the equation
$P=(p_0s^2+p_1st+p_2t^2)xy$. This leads to the following
six equations:
$$
\matrix{
\pi_1={}&    &   & \hfill{}-p_2(\xi_1\eta_2+\xi_2\eta_1) \cr
\pi_2={}&p_0\xi_1\eta_1 \hfill &      & \hfill{}-p_2\xi_2\eta_2 \cr
\pi_3={}&p_0(\xi_1\eta_2+\xi_2\eta_1) \hfill& {}+p_1\xi_2\eta_2 \hfill\cr
\pi_4={}&p_0(\xi_1\eta_3+\xi_2\eta_2)+\xi_3\eta_1  \hfill
                   &{}+p_1(\xi_2\eta_3+\xi_3\eta_2) \hfill
                                        & {}+p_2\xi_3\eta_3  \hfill\cr
\pi_5={}&p_0(\xi_2\eta_3+\xi_3\eta_2)+\xi_4\eta_1 \hfill
    &{}+p_1(\xi_3\eta_3+\xi_4\eta_2) \hfill&{}+p_2 \xi_4\eta_3  \hfill\cr
\pi_6={}&p_0(\xi_3\eta_3+\xi_4\eta_2) \hfill  &{}+p_1\xi_4\eta_3  \hfill &  \cr
}
$$
If we write a matrix with the coefficients of the $p_i$ in the columns
with rows coming from the equations $\pi_m$ we find that the first 
$k+1$ rows form a skew symmetric matrix. This is due to the specific
choices made in the above proof. One can also get any other  block to 
be skew symmetric by using the lifting conditions. In this example they are
$p_0\eta_1+p_1\eta_2+p_2\eta_3=0$, $p_0\xi_1+p_1\xi_2+p_2\xi_3=0$
and $p_0\xi_2+p_1\xi_3+p_2\xi_4=0$.
From the skew symmetry we can conclude:

\proclaim Proposition \num.
If $e_y\leq e_x<b$ then the $b-1$ equations $\pi_m$ coming from
the equation $P=(\sum_{j=0}^k p_js^{k-j}t^k)xy$, where $b+k=e_x+e_y$,
satisfy $b-k-1$ linear relations
$\sum_{j=0}^k  p_j\pi_{i+j} =0$, for $0< i < b-k$. 

\roep Example \num: case II. 
Let $b=4$, $e_x=5$ and $e_y=3$. Consider the equation
$P=(p_0s^4+p_1s^3t+p_2s^2t^2+p_3st^3+p_4t^4)xy$. This leads to the following
three equations:
$$
\matrix{
\pi_1=\rho_0\xi_1+\rho_1\xi_2 &
{}-p_2\xi_2\eta_1-p_3(\xi_2\eta_2+\xi_3\eta_1)
     p_4(\xi_3\eta_2+\xi_4\eta_1) \cr
\pi_2=\rho_0\xi_2+\rho_1\xi_3&
  {}+ p_0\xi_1\eta_1+p_1\xi_2\eta_1 -p_3\xi_3\eta_2-p_4\xi_4\eta_2\cr
\pi_3=\rho_0\xi_3+\rho_1\xi_4& 
  +{}p_0(\xi_1\eta_2+\xi_2\eta_1)+p_1(\xi_2\eta_2+\xi_3\eta_1) 
 +p_2\xi_3\eta_2 \cr
}
$$
\endroep
 

\roep {\rm(\num)} Hyperelliptic cones {\rm(cf.~\cite{\rfsth})}.
Let $X$ be the cone over a hyperelliptic curve $C$ embedded with a
line bundle $L$ of degree $d\geq 2g+3$. Then $\dim T^1_X(-1)=2g+2$.
The curve lies on a scroll of degree $d-g-1$ as curve of type
$2H-(d-2g-2)R$. The number of rolling factors equations
is $d-2g-3$, so we have at least as many equations as variables
if $d>4g+4$. In that case only conical deformations exist, so 
all deformations in negative degree are obstructed.

The easiest case to  describe is $L=ng^1_2$. 
The curve $C$ has an affine equation
$y^2=\sum_{k=0}^{2g+2} p_kt^k$, which gives the bihomogeneous
equation $(\sum_{k=0}^{2g+2} p_ks^{2g+2-k}t^k)x^2-y^2=0$.
The line bundle $L$ embeds $C$ in a scroll $S(n,n-g-1)$,
and there are $2n-2g-1$ rolling factors equations $P_m$,
coming from $p(s,t)x^2-y^2$.
The lifting matrix is a block diagonal matrix with the $y$-block
equal to $-2I_{n-g-2}$, and the $x$-block a $(n-2g-3)\times(n-1)$
matrix, so the dimension of the space of
lifting deformations of the scroll is $2g+2$ if $n\geq 2g+3$.
If $n\leq 2g+3$, the $x$-block is not present, and all $n-1$
$\xi$-deformations lift. Furthermore there are $2g+3-n$ pure
rolling factors deformations. This shows again that
$\dim T^1_X(-1)=2g+2$.

\proclaim Proposition \num.
If $n\geq 2g+3$ the base space in negative degrees is
a zero-dimensional complete intersection of $2g+2$
quadratic equations.

\roep Proof.
We may assume that the highest coefficient $p_{2g+2}$ in
$p(s,t)$ equals $1$. The lifting equations  allow now to
eliminate the variables $\xi_{2g+3},\dots,\xi_{n-1}$.
The base equations $\pi_m$ involve only the $\xi_i$ and are
therefore not linearly independent. Because $p_{2g+2}=1$
we can discard all $\pi_m$ with $m>2g+2$. The first
$2g+2$ equations involve only the first $2g+2$ variables.
This shows that we have the same system of equations for 
all $n\geq 2g+3$. As we know that there are no deformations over a
positive dimensional base, we conclude that the base space
is a complete intersection of $2g+2$ equations.
\qed

\roep Remark \num.
The fact that the system of equations above defines a complete intersection 
can also be seen directly. In fact we have the following result:

\proclaim Lemma \num. The system of $e=b-1$ equations $\pi_m$ 
in $e-1$ variables $\xi_i$ coming from
one polynomial $P_{b-2}(s,t)x^2$  is a zero-dimensional complete 
intersection if and only if $P_{b-2}(s,t)$ has no multiple roots.

\edef\labprop{{\the\secno.\the\thmno}}
\roep Proof.
First we note that there are only $b-2$ linearly independent
equations.
We put $\xi_i=s^{b-i-1}t^i$.
Then 
$$
\displaylines{\qquad
  \pi_m=\sum_{k=0}^{b-2}(m-k-1)p_ks^{2b-k-m-2}t^{k+m-2}
\hfill\cr\hfill{}
    =s^{b-m}t^{m-2}\bigl(\sum(b-2-k)p_ks^{b-2-k}t^k + 
    (m+1-b)\sum p_ks^{b-2-k}t^k\bigr)\qquad\;.}
$$
The form $P(s,t)$ has multiple roots if and only if  
$P(s,t)$ and $s{\partial \over \partial s}P(s,t)$ have a common zero
$(s_0\cn t_0)$. Then $\xi_i=s_0^{b-i-1}t_0^i$ is a nontrivial
solution to the system of equations.

We show the converse by induction. One first checks that a linear 
transformation in $(s \cn t)$
does not change the isomorphism type of the ideal. We apply a
transformation such that $s=0$ is a single root of $P$, so 
$p_0=0$ but $p_1\neq 0$. The equations $\pi_2$, \dots, $\pi_{b-1}$
now do  not involve the variable $\xi_1$ and are by the induction
hypotheses a complete intersection in $e-2$ variables, so their zero
set is the $\xi_1$-axis with multiple structure. The equation
$\pi_1$ has the form $-p_1\xi_1^2+\dots$, so the whole system has
a zero-dimensional solution set.
\qed

\comment
In the general case of a hyperelliptic cone we can write $L=ng^1_2+D$ with 
$D$ an effective divisor of lowest possible degree: $\deg D=g+1-e$.
If we describe $C$ with the affine equation 
$y^2=\sum p_kt^k$, then the point set $D$ can be given by
the ideal $(U(t),y-V(t))$, where the polynomial $U(t)$ has as
zeroes the $t$-coordinates of the points in $D$. One can write 
$P=V^2+UW$ for some polynomial $W$. Sections of $L$ can be
represented by the forms $x_i=Ut^i$, $0\leq i\leq n$
and $y_i=(y+V)t^i$, $0\leq i \leq n-e$.
We can write $P-y^2=UW+2V(y+V)-(y+V)^2$ and this gives the
bihomogeneous equation
$$
W_{g+1+e}(s,t)x^2+2V_{g+1}(s,t)xy-U_{g+1-e}(s,t)y^2
$$
on the scroll $S(n,n-e)$. There are no pure rolling factors deformations
if $n>g+1+e$. If $\deg L>4g+4$ the base space is zero-dimensional,
but the equations involve terms coming from $U$, $V$ and $W$.
The terms from $W$ give only $g+1+e$ linearly independent expressions,
those from $V$ and $U$ respectively $g+1$ and $g+1-e$. Therefore
each base equation can be expressed as linear combination of
$3g+3$ expressions.
\endcomment
\roep Remark \num.
For $\deg L= 4g+4$ the base space is a cone over $2^{2g+1}$ points
in a very special position: there exist $2g+2$ hyperplanes $\{l_i=0\}$
such that the base is given by $l_i^2=l_j^2$ \cite{\rfsth}. We can make this
more explicit in the case $L=(2g+2)g_2^1$. 
Again the $y$-block of the lifting matrix is a multiple of the
identity, but now there is also one rolling factors deformation 
parameter $\rho$. More generally, we look the equations
coming from $p(s,t)x^2$ with $\deg p = b =e$. We get base
equations $\Pi_m=\rho \xi_m+\pi_m$, where $\pi_m$ is a 
quadratic equation in the $\xi$-variables only.
One solution is clearly $\xi_i=0$ for all $i$. To find the others
we eliminate $\rho$:
$$
\Rank \pmatrix {
\pi _1 & \pi_2 & \dots & \pi_{e-1} \cr
\xi _1 & \xi_2 & \dots & \xi_{e-1}}
\leq 1\;. \eqno{(**)}
$$
The equations $\Pi_m$ can be changed by changing $\rho$,
but this system is independent of such changes.
Write inhomogeneously $p(t)=p_0+p_1t+\dots+p_{e-1}t^{e-1}+t^e=
\prod (t-\al_i)$, where the $\al_i$ are the roots of $p(t)$.

\proclaim Lemma \num.
The $e$ points $P_i=(1\cn \al_i \cn \al_i^2 \cn \cdots \cn \al_i^{e-2})$
are solutions to the system $(**)$.

\roep Proof.
Let $\al$ be a root of $p$ and insert $\xi_i=\al^{i-1}$ in the
system $(**)$. We simplify the matrix by 
column operations: subtract  $\al$ times the  $j$th column from
the $(j+1)$st column, starting at the end.
The matrix has clearly rank $1$, if $\pi_{j+1}(\al)-\al\pi_j(\al)=0$,
where $\pi_j(\al)$ is the result of substituting $\xi_i=\al^{i-1}$
in the equation $\pi_j$. The coefficient $p_k$ occurs in 
$\pi_j(\al)$ in the term $lp_k\al^{j+k-2}$ for some integer $l$, and in 
the term $(l+1)p_k\al^{j+k-1}$ in $\pi_{j+1}(\al)$.
Therefore  $\pi_{j+1}(\al)-\al\pi_j(\al)=-\sum p_k \al^{j+k-1}
=-\al^{j-1}p(\al)=0$.
\qed

The remaining solutions are found in the following way.
Divide the set of roots into two subsets $I$ and $J$.
The points $P_i$ lie on a rational normal curve. Therefore the
points $P_i$ with $i\in I$ span a linear subspace $L_I$ of dimension 
$|I|-1$. 

\proclaim Claim. The intersection point $P_I:=L_I\cap L_J$ is a solution
to $(**)$.

The proof  is a similar but more complicated computation.
We determine here only the point $P_I$. 
The condition that the point $\sum_{i\in I} \la_i P_i$ lies in $L_J$
is that
$$
\Rank \pmatrix {
\sum \la_i & \dots & \sum \la_i\al_i^{e-2} \cr
 1 & \dots & \al_{j_1}^{e-2} \cr
\vdots & & \vdots \cr
 1 & \dots & \al_{j_{|J|}}^{e-2} }
= |J|\;.
$$
We find the resulting linear equations on the $\la_i$
by extending the matrix to
a square matrix by adding $|I|-2$ rows of points on the rational normal
curve, for which we take roots. Then only two $\la_i$ survive, and
they come with a Vandermonde determinant as coefficient.
Upon dividing by common factors we get 
$(\prod_{i\neq i_1,i_2}(\al_{i_1}-\al_i))\la_{i_1} 
+(\prod_{i\neq i_1,i_2}(\al_{i_2}-\al_i))\la_{i_2}=0 $. We multiply
with $\al_{i_1}-\al_{i_2}$. Noting that
$\prod_{i\neq i_1}(\al_{i_1}-\al_i)=p'(\al_{i_1})$ (with $p'(t)$ the 
derivative of $p(t)$) we
get $p'(\al_{i_1})\la_{i_1}=p'(\al_{i_2})\la_{i_2}$.

\comment
\roep Proof.
We first compute $P_I$.
The condition that the point $\sum_{i\in I} \la_i P_i$ lies in $L_J$
is that
$$
\Rank \pmatrix {
\sum \la_i & \dots & \sum \la_i\al_i^{e-2} \cr
 1 & \dots & \al_{j_1}^{e-2} \cr
\vdots & & \vdots \cr
 1 & \dots & \al_{j_{|J|}}^{e-2} }
= |J|\;.
$$
This gives a system of linear equations 
on the $\la_i$. We can find these equations by extending the matrix to
a square matrix by adding $|I|-2$ rows of points on the rational normal
curve, for which we take roots. Then only two $\la_i$ survive, and
they come with a Vandermonde determinant as coefficient.
Upon dividing by common factors we get 
$(\prod_{i\neq i_1,i_2}(\al_{i_1}-\al_i))\la_{i_1} 
+(\prod_{i\neq i_1,i_2}(\al_{i_2}-\al_i))\la_{i_2}=0 $. We multiply
with $\al_{i_1}-\al_{i_2}$. Noting that
$\prod_{i\neq i_1}(\al_{i_1}-\al_i)=p'(\al_{i_1})$ we
get $p'(\al_{i_1})\la_{i_1}=p'(\al_{i_2})\la_{i_2}$.

We insert the coordinates of the point $P_I$ into the system $(**)$.
We first suppose that $|I|=2$, so we have two roots $\al$, $\be$
and the point $P=\la P_{\al}+\mu P_{\be}$ with $p'(\al)\la=p'(\be)\mu$.
We simplify the matrix with column operations, again subtracting  
$\al$ times the  $i$th column from the $(i+1)$st column, and repeating
with $\beta$ on the resulting matrix. The first row
consists now of zeroes, except in the first two columns, while
the second row has entries $ \pi_{i+2}-(\al+\be)\pi_{i+1}+\al\be\pi_i$.
This is a quadratic polynomial in $\la$ and $\mu$.
By writing it as $ (\pi_{i+2}-\al\pi_{i+1})-\be(\pi_{i+1}-\al\pi_i)$
and using the computation of the previous Lemma
we recognise that the coefficient of $\la^2$ vanishes.
The computation that also the coefficient of $\la\mu$ vanishes,
is left to the reader. 
To prove that we have a solution it therefore suffices to show 
that the first minor vanishes.
Note that 
$$
\bmatrix{
\la+\mu & \la\al+\mu\be \cr
\pi_1 & \pi_2 }
=
\bmatrix{
\la & \mu \cr
\be\pi_1-\pi_2 & \pi_2-\al \pi_1 }
\;.
$$
We write $\pi_i=A_i\la^2+B_i\la\mu+C_i\mu^2$.
We determine the terms with $\la^2\mu$ and $\la\mu^2$.
They are
$$
\la\mu((A_2+B_2-\al B_1-\be A_1)\la+(C_2+B_2-\be B_1-\al C_1)\mu) \;.
$$
As we know that $\al A_1-A_2=0$ it remains to see that
$-2A_2+B_2-\al B_1-\be A_1+3\al A_1=(\be-\al)p'(\al)$.
We look at the monomial $p_kt^k$. We may choose the equations $\pi_m$
such that we  get for $k\geq 3$ in $\pi_1$ the contribution
$-\xi_2\xi_{k-1}-\dots - \xi_{k-1}\xi_2$ and
$-\xi_3\xi_{k-1}-\dots - \xi_{k-1}\xi_3$  in $\pi_2$.
Then $A_1=-(k-2)\al^{k-1}$,  $A_2=-(k-3)\al^{k}$,
$B_1=-2(\al\be^{k-2}+\dots+\al^{k-2}\be)$ and
$B_2=-2(\al^2\be^{k-2}+\dots+\al^{k-2}\be^2)$.
This gives 
$-2A_2+B_2-\al B_1-\be A_1+3\al A_1=
(\be-\al)k\al^{k-1}$.
For $p_0$ we get only a contribution from $\pi_2$: $A_2=1$, $B_2=2$.
For the remaining small $k$ we may assume that we have 
$-\xi_1\xi_{k}-\dots - \xi_{k}\xi_1$, 
resp.~$-\xi_2\xi_{k}-\dots - \xi_{k}\xi_2$.
This makes
$A_1=-k\al^{k-1}$,  $A_2=-(k-1)\al^{k}$,
$B_1=-2(\be^{k-1}+\al\be^{k-2}+\dots+\al^{k-2}\be+\al^{k-1})$ and
$B_2=-2(\al\beta^{k-1}+\dots+\al^{k-1}\be)$.
This gives also
$-2A_2+B_2-\al B_1-\be A_1+3\al A_1=
(\be-\al)k\al^{k-1}$.

If we have more than two roots, the minors are cubic polynomials
in the $\la_i$ and we know that all coefficients vanish except
possibly those of $\la_i\la_j\la_k$ with pairwise different
indices.
So it suffices to check the case of three roots $\al$, $\be$ and $\ga$.
By suitable operations on the matrix the first three columns give 
after transposing
$$
\pmatrix{
\la &(\be-\ga)(\pi_3-(\be+\ga)\pi_2+\be\ga\pi_1) \cr
\mu &(\ga-\al)(\pi_3-(\ga+\al)\pi_2+\ga\al\pi_1) \cr
\nu &(\al-\be)(\pi_3-(\al+\be)\pi_2+\al\be\pi_1) 
}\;.
$$
As there is no $\la\nu$-term in $\pi_3-(\ga+\al)\pi_2+\ga\al\pi_1$
and no $\la\mu$-term in $\pi_3-(\al+\be)\pi_2+\al\be\pi_1$ we are 
reduced to check that the $\la\mu\nu$-term vanishes
in $(\mu(\al-\be)+\nu(\ga-\al))(\pi_2-\al\pi_1)$, which is similar to the
computation above.
\qed
\endcomment

We write out the equations for $e=5$:
$$
\displaylines{
\rho\xi_1- p_1\xi_1^2-2p_2\xi_1\xi_2-p_3\xi_2^2-2p_4\xi_2\xi_3 
                         -p_5(2\xi_2\xi_4+\xi_3^2) \cr
\rho\xi_2+ p_0\xi_1^2-p_2\xi_2^2-p_4\xi_3^2-2p_5\xi_3\xi_4 \cr
\rho\xi_3+ 2p_0\xi_1\xi_2+p_1\xi_2^2-p_3\xi_3^2-p_5\xi_4^2 \cr
\rho\xi_4+ p_0(2\xi_1\xi_3+\xi_2^2)+2p_1\xi_2\xi_3+p_2\xi_3^2 
                         +2p_3\xi_3\xi_4+p_4\xi_4^2 
}
$$
Let $\al$ be a root of $p_0+p_1t+p_2t^2+p_3t^3+p_4t^4+t^5$,
and $\be$, \dots, $\ep$ the remaining roots.
Write $\sigma'_i$ for the $i$th symmetric function of these four roots.
Then a solution is $\xi_i=\al^{i-1}$, $\rho=\al^4-\al^3\sigma'_1
-\al^2\sigma'_2-\al\sigma'_3+\sigma'_4$.
Given two roots $\al$ and $\be$ we get a solution
$\xi_i=(\ga-\be)(\de-\be)(\ep-\be)\al^i+(\al-\ga)(\al-\de)(\al-\ep)\be^i$.
To write $\rho$ we set $\mu=(\ga-\be)(\de-\be)(\ep-\be)$,
$\la=(\al-\ga)(\al-\de)(\al-\ep)$ and $\sigma_i''$ the $i$th symmetric
function in $\ga$, $\de$ and $\ep$. Then 
and $\rho=\mu
(\al^4-\al^2(\al+2\be)\sigma_1''-\al^2\sigma_2''-\al\sigma_3'')
+\la
(\be^4-\be^2(\al+2\al)\sigma_1''-\be^2\sigma_2''-\be\sigma_3'')
$.
The hyperplane through $(1\cn0\cn0\cn0\cn0)$, $P_\ga$, $P_\de$ and 
$P_\ep$ is 
$l^-_{\al\be}=\sigma_3''\xi_1-\sigma_2''\xi_2+\sigma_1''\xi_3-\xi_4$.
In it lie also $P_{\ga\de}$, $P_{\ga\ep}$, $P_{\de\ep}$ and $P_{\al\be}$.
The hyperplane containing the remaining points is 
$l^+_{\al\be}=\rho-(\al+\be)l^-_{\al\be}+2\sigma_3''\xi_2+2\al\be\xi_3$.
We put $l_a=\rho-2\sigma'_4\xi_1+2\sigma'_3\xi_2+2\al\sigma'_1\xi_3
 -2\al\xi_4$.
Then $l_\al^2-l_\be^2=4(\al-\be)l^-_{\al\be}l^+_{\al\be}$.

\section Tetragonal curves.
 
An extension of a  canonical curve 
yields a surface with the given canonical curve
as hyperplane section. Surfaces with canonical hyperplane sections
were studied in Dick Epema's thesis \cite{\rfep}.
Only a limited list of surfaces can occur.

\proclaim Theorem \num {\rm (\cite{\rfep}, Cor.~I.5.5 and Cor.~II.3.3)}.
Let $W$ be a surface with canonical hyperplane sections. Then one
of the following holds:
\item {\rm (a)} $W$ is a $K3$ surface with at most rational double points
          as singularities,
\item {\rm (b)} $W$ is a rational surface with one minimally elliptic
          singularity and possibly rational double points,
\item {\rm (c)} $W$ is a birationally ruled surface over an 
          elliptic curve $\Gamma$
          with as non-rational singularities either
\itemitem {\rm i)} two simple elliptic singularities with exceptional divisor
             isomorphic to $\Gamma$, or
\itemitem {\rm ii)} one Gorenstein singularity with $p_g=2$,
\item {\rm (d)} $W$ is a birationally ruled surface over a curve $\Gamma$
          of genus $q\geq2$
          with one non-rational singularity with $p_g=q+1$, whose
          exceptional divisor contains exactly one non-rational curve
          isomorphic to $\Gamma$.

Case (c) occurs for bi-elliptic curves (see below). 
If we exclude them and curves of low genus on Del Pezzo surfaces,
then all extensions of tetragonal curves
are of rolling factors type. The surface $W$ has
therefore to occur in our classification of complete intersection
surfaces on scrolls. In particular, $K3$ surfaces can only occur
if $b_1\leq b_2+4$.
This has consequences for deformations of tetragonal cones.

\proclaim Proposition \num.
Pure rolling factors deformations are always unobstructed.
If $e_3>0$ and $b_1>b_2+4$ the remaining deformations are obstructed.

\roep Proof.
The first statement follows directly from the form of the equations.
For the second we note that the total space of
a nontrivial one-parameter deformation
of  a scroll with $e_3>0$ is a scroll with $e_4>0$.
\qed

By taking hyperplane sections of a general element in each of the families
of the classification we obtain for all $g$ tetragonal curves with
$b_1\leq b_2+4$  lying on $K3$ surfaces (with at most rational
double points). To realise the other types of surfaces
we give a construction, which goes back to 
\cite{\rfduv}.  His construction was generalised to the 
non-rational case in \cite{\rfep}. In our
situation  we want a given curve to be a hyperplane section.
A general construction for given hyperplane sections of 
regular surfaces is given in \cite{\rfwacy}.

\proclaim Construction \num.
Let $Y$ be a surface containing the curve $C$ and let $D\in
|-K_Y|$ be an anticanonical divisor. 
Let $\wt Y$ be the blow up of $Y$ in the scheme $Z=C\cap D$.
If  the linear subsystem 
$\cc'$ of $|C|$ with base scheme $Z$ has dimension $g$, 
it associated map contracts $D$ and blows down $\wt Y$ to
a surface $\wl Y$ with $C$ as canonical hyperplane section.

Let $\ciz$ be the ideal sheaf of $Z$. Then we have
the exact sequence
$$
0 \lra \sier Y
\lra \ciz  \sier Y(C)
\lra \sier C(C-Z) \lra 0
$$
and by the adjunction formula $\sier C(C-Z)=K_C$.
If $ h^0(\ciz  \sier Y(C))=g+1$ then the map
$H^0(\ciz  \sier Y(C))\lra H^0(K_C)$ is surjective, a condition
which is automatically satisfied if $Y$ is a regular surface.
This yields that the special hyperplane section is the curve $C$
in its canonical embedding.

Suppose that $Y$ is not regular. By Epema's classification
$Y$ is then a birationally ruled surface, over a curve
$\Gamma$ of genus $q$.
Let $\wt C$ be the strict transform of $C$ on $\wt Y$ and $\wl C$ its
image on $\wl Y$. Then $H^0(\ciz  \sier Y(C))=H^0(\sier {\wt Y}(\wt C))
=H^0(\sier {\wl Y}(\wl C))$.
We look at the exact sequence 
$$
0 \lra \sier {\wl Y}
\lra \sier {\wl Y}(\wl C)
\lra \sier {\wl C}(\wl C)=K_C\lra 0 \;.
$$
We compute $H^1( \sier {\wl Y})$ with the spectral
sequence for the map $\pi\colon \wt Y \lra \wl Y$.
This gives us the long exact sequence
$$
0 \lra H^1( \sier {\wl Y})
\lra H^1( \sier {\wt Y}) 
\lra H^0( R^1\pi_*\sier {\wt Y})
\lra H^2( \sier {\wl Y})
\lra 0
$$
in which $\dim H^1( \sier {\wt Y}) =q$.
We choose $D$ in such a way that the composed map
$H^1(\sier \Gamma)\lra R^1\pi_*\sier {\wt Y} \lra
H^1(\sier{\wt D})$, where $\wt D$ is the exceptional divisor
of the map $\pi$, is injective.
Then the map $H^0(\ciz  \sier Y(C))\lra H^0(K_C)$
is surjective.

To apply the construction we need a surface on which the curve
$C$ lies. In the tetragonal case
a natural candidate is the surface $Y$ of type
$2H-b_1R$ on the scroll. 

We first assume that $e_1<b_1$, so there are no pure rolling factors
deformations coming from the first equation on the scroll.
The canonical divisor of the scroll $S$ is $-3H+(b_1+b_2)R$
\cite{\rfschr, 1.7}. So an anticanonical divisor on $Y$ is of type
$H-b_2R$.
Let $T=\tau_{e_1-b_2}(s,t) x +\tau_{e_2-b_2}(s,t)y+\tau_{e_3-b_2}(s,t) z$
be the equation of such a divisor.
Sections of $ \ciz  \sier Y(C)$ are $Q$ (which defines $C$), and
$x_iT=s^{e_1-i}t^ixT$, $y_iT$ and $z_iT$. With coordinates $(t\cn x_i\cn y_i
\cn z_i )$ on $\P^{g}$ we get by rolling factors 
$b_2+1$ equations $\wt Q_m$ from the relation 
$Q(\tau_{e_1-b_2}(s,t) x +\tau_{e_2-b_2}(s,t)y+\tau_{e_3-b_2}(s,t) z)
=(Q_{1,1}x^2+\cdots+Q_{3,3}z^2)T$. As $t$ is also a coordinate on the
four-dimensional scroll, which is the cone over $S$, we can write the 
equation on the scroll as
$$
Q_{1,1}x^2+\cdots+Q_{3,3}z^2-
(\tau_{e_1-b_2} x +\cdots+\tau_{e_3-b_2} z)t\;.
$$

We analyse the resulting singularities.
If $Y$ is a rational surface, we have a anticanonical divisor $D$
which has arithmetic genus 1, giving a minimally elliptic singularity
on the total space of the deformation. 

If $Y$ is a ruled surface over a hyperelliptic curve $\Gamma$,
then $D$ passes through the double locus. This gives an exceptional
divisor with $\Gamma$ as only non-rational curve.

\roep Example \num.
Let  $(e_1,e_2,e_3;b_1,b_2)=(3n-2,2n-2,n-2;4n-4,2n-4)$.
If the coefficient of $xz$ does not vanish, we may bring the equation $P$
onto the form $xz-y^2$.
The second equation has the form $z^2+q_nzy+q_{2n}zx+q_{3n}xy+q_{4n}x^2$
from which $z$ may be eliminated to obtain a quartic equation for $y$.
The case of a cyclic curve $y^4+q_{4n}x^4$ is a special instance.
The  equation $P$ gives  a square lifting matrix
in which the antidiagonal blocks are square unit matrices.
Therefore the only deformations are pure rolling factors deformations,
coming from the second equation, in number $(n+3)+3={g+3\over 6}+6$.
We have $T=\tau_{n+2}x+\tau_2y$. The section $(0\cn0\cn1)$ is always
a component of $D$. If $t_2\not\equiv0$ we have a cusp singularity,
but if $\tau_2\equiv0$ the section occurs with multiplicity 2 in $D$.

If however the coefficient of $xz$ vanishes, the surface $Y$ is singular.
After a coordinate
transformation its equation is $y^2+p_{2n}x^2$, the other equation being
$z^2+q_{3n}xy+q_{4n}x^2$. In this case the lifting matrix has 
(up to a factor $\frac12$) the 
following block structure
$$
\pmatrix
{ \Pi & 0 & 0 \cr
   0  & I & 0 \cr
   0  & 0 & 0 \cr
   0  & 0 & I}
$$
so there are $2n$ $\xi$-deformations, on which we have $4n-5$ base
equations coming from the equation $P$. Of these are only $2n$ linearly
independent, defining a zero-dimensional complete intersection
(see Lemma~\labprop). These deformations are therefore obstructed,
leaving us again with only the pure rolling factors deformations.
The curve $D$ consists of the double locus and in general
$2n+4$ lines.
\endroep

The same computation as above works for bielliptic cones. In that case
one has a deformation of weight $-2$. The total space
is a surface in weighted projective space $\P(1,\dots,1,2)$. Replacing the
deformation parameter $t$ by $t^2$ we get a surface in ordinary $\P^g$.
This is a surface with two simple elliptic singularities.
The most general surface of this type is the intersection
of our elliptic cone with one dimensional vertex with the hypersurface
given by
$$
\wt Q= z^2+Q(x_i,y_i)+tl(x_i,y_i)+at^2\;,
$$
where $l(x_i,y_i)$ is a linear form in the coordinates $x_i$, $y_i$.
If the coefficient $a$ vanishes, we get a surface with one singularity
with $p_g=2$.
The construction above gives an equation of the form
$\wt Q= z^2+\cdots+azt$, which after a coordinate transformation
becomes $z^2+\cdots-\frac14a^2t^2$.

\proclaim Proposition \num.
For bielliptic cones of genus $g>10$ the only deformations of negative weight
are pure rolling factors deformations.

\roep Proof.
Each infinitesimal deformation of the bielliptic cone induces
an infinitesimal deformation of the cone over the projective cone 
over the elliptic curve. The same holds therefore for complete
deformations of negative weight.
It is well-known that the cone over an elliptic curve of degree at least 10
has only obstructed deformations of negative weight. Therefore the
deformation of the elliptic cone is trivial and the
only possibility is to deform the last quadratic equation.
\qed

On the other hand, non-scrollar extension do occur for bielliptic 
curves with
$g\leq 10$ and for tetragonal curves on Del Pezzo surfaces.

\roep Example \num.
A bielliptic curve of genus 10 lies on the projective cone over an
elliptic curve of degree $9$. Such a cone is can be smoothed to the
triple Veronese embedding of $\P^2$.  
Let $W$ be a $K3$ surface of degree 2, a double cover of $\P^2$ branched
along a sextic curve. We re-embed $W$ with $|3L|$, where $L$ is the pull-back
of  a line on $\P^2$. The image lies on the cone over the
Verones embedding. A hyperplane section through the vertex of the cone is
a bielliptic curve, whereas the general hyperplane section has
a $g^2_6$.  This example, due to \cite{\rfdonmor}, is the only case
where the gonality of smooth curves in a base-point-free ample linear
system on a $K3$ surface is not constant \cite{\rfcp}.

\endroep

Now we look at the case that also the first set of equations
admit pure rolling factors deformations.

\proclaim Lemma \num. If $e_1\geq b_1$ then $e_1\leq b_1+2$ and
$b_1\leq b_2+4$.

\roep Proof.
Under the assumption $e_1\geq b_1$ we have 
$e_2+e_3-2\leq b_2\leq 2e_3$ so $e_2\leq e_3+2$ and
$b_1\leq 2e_2\leq e_2+e_3+2\leq b_2+4$. Furthermore
$e-1=b_1+b_2+2-e_2-e_3\leq b_1+2$.
\qed

It is now easy to list all 18 possibilities, ranging
from $(2e+2,e,e;2e,2e)$ to $(2e+4,e+2,e;2e+4,2e)$.
A look at the table of tetragonal $K3$ surfaces
reveals that all possibilities are
realisable as  special sections of $K3$-surfaces;
e.g., the hyperplane section $x_{e+2}=y_0$
of a $K3$ with invariants $(e+2,e+2,e+2,e;2e+4,2e)$ yields
the last case.

On the other hand, every family of $K3$ surfaces contains degenerate
elements with singularities of higher genus. Those can be constructed
with Epema's construction and in fact he gives rather complete
results for quartic hypersurfaces \cite{\rfep}. The classification
of such surfaces is  due to \cite{\rfrohn} and is quite involved.
In those cases the rational or ruled surfaces on which the 
canonical curve lies are not evident.
For pure rolling factors extensions the situation is better;
in fact, we can make the following simple observation.

\proclaim Proposition \num.
Let $W$ be a pure rolling factors extension of tetragonal curve,
which is not bi-elliptic. It lies on the cone over the
3-dimensional scroll $S$ with vertex in $p=(0\cn \dots \cn 0\cn 1)$
and the projection from the point $p$ yields a surface $Y\subset
\P^{g-1}$ on which $C$ lies.

\roep Example \num. 
If $b_1>e_1$ then $X$ lies on the cone over the surface $Y$
on the scroll and the projection is just this surface $Y$, so
we get the construction described above.

\roep Example \num. Consider the curve with invariants $(8,4,2;8,4)$.
In general a pure rolling factors extension leads to a 
$K3$-surface (with an ordinary double point). 
It is the case $e=3$ of $(e+5,e+1,e-1,e-3;2e+2,2e-1)$
from the table; the singularity appears because the section
$(0\cn 0\cn 0\cn1)$ is contracted.
To find the equation of $Y$ on the scroll we  have to eliminate the
last coordinate $w$. The deformed equation $\wt P$ is  $P+axw$, while
$\wt Q=Q+byw+c_2(s,t)xw$ with $a$ and $b$ nonzero constants. The equation
of $Y$ is therefore $(by+c_2(s,t)x)P-axQ$, which defines a
divisor of type $3H-(b_1+b_2)R$ on the scroll.

\roep Example \num: the case $(2e+2,e,e;2e,2e)$.
We first derive a normal form for the equations $P$ and $Q$.
We start with the restriction to $x=0$. We have a pencil of quadrics
so we may choose the first equation as $y^2$ and the second as $z^2$.
We get:
$$
\eqalign
{P\colon\;& y^2+p_{e+2}xz+p_{2e+4}x^2 \cr
Q\colon\;& z^2+q_{e+2}xy+q_{2e+4}x^2 \;.\cr}
$$
There are $3+3$ pure rolling factors deformations:
$$
\eqalign
{\wt P\colon\;& P+(\rho_0 s^2+\rho_1 st +\rho_2t^2)x \cr
\wt Q\colon\;& Q+(\tau_0 s^2+\tau_1 st +\tau_2t^2)x \;.\cr}
$$
If the polynomials $\rho:=\rho_0 s^2+\rho_1 st +\rho_2t^2$ and
$\tau:=\tau_0 s^2+\tau_1 st +\tau_2t^2$ are proportional, so
$\la \rho + \mu \tau =0$, then the surface $Y$ is the surface 
$\la P + \mu Q=0$ from the pencil. In general the anticanonical
divisor $D$ contains the two sections given by $x=0$,
$\la y^2+\mu z^2=0$ and the singularity on the deformation
is a cusp singularity. If $\rho$ and $\tau$ have $0\leq\gamma<2$
roots in common, the projected surface is a divisor of type
$2H - (2e-2+\gamma)R$. In general we get a simple elliptic singularity.

To describe the remaining deformations
we look at the lifting matrix, which is a block matrix
$$
\pmatrix{
0 & 2I & 0 \cr
\Pi_{e+2} & 0 & 0 \cr
\Xi_{e+2} & 0 & 0 \cr
0 & 0 & 2I \cr}
$$
of size $(4e-4)\times (4e-1)$. Its rank is $2e-2$ if $p_{e+2}$ and $q_{e+2}$
both vanish identically, and lies between $3e-3$ and $4e-4$
otherwise. The solution space has dimension $\gamma\geq 3$
with strict inequality iff the polynomials $p_{e+2}$ and $q_{e+2}$ have 
$\gamma$ roots in common.
The $\eta$ and $\zeta$ deformations vanish. Therefore the base
equations depend only on $p_{2e+4}$ and $q_{2e+4}$. They are
$2(2e-1)$ quadratic equations on $2e+1+6$ variables, which may or may not
have solutions.
\endroep

We now turn to the other deformations in general.
A dimension count shows that the general tetragonal curve
of genus $g>15$ cannot lie on a $K3$ surface, so the deformations
are obstructed. 
For a general tetragonal cone we have that $\dim T^1_X=9$. There are
$(b_1-1)+(b_2-1)=g-7$ quadratic base equations. 
Compare this with the dimension of $T^2$:

\proclaim Theorem \num {\rm(\cite{Wahl, Thm.~5.9})}.
Let $X$ be a tetragonal cone with $e_3>0$. Then 
$\dim T^2_X(-k)=0$ for $k>2$ and 
$\dim T^2_X(-k)=g-7$ if $b_2>0$. If $b_2=0$, then 
$\dim T^2_X(-k)=2(g-6)$.

In particular, if $g>15$ we have more equations than variables and in general
there are no solutions. For special moduli solutions do exist and one
expects in general exactly one solution.

\roep {\rm(\num)} The case $g=15$.
Consider the most general situation, 
of equal invariants: $e_1=e_2=e_3=4$, $b_1=b_2=5$.
In this case there are no pure rolling factor deformations and no
lifting conditions. 

\proclaim Proposition \num. The general tetragonal curve
with $e_1=e_2=e_3=4$, $b_1=b_2=5$ is hyperplane section of 
$256$ different $K3$ surfaces.

\roep Proof. We have $8$ homogeneous quadratic equations in $9$ variables,
which define a complete intersection of degree $2^8$. 
We give an explicit example. Take the curve, given by
the equations 
$$
\displaylines{
(s^3+t^3)x^2+(s^3+2t^3)y^2+(s^3-2t^3)z^2  \cr
(s^2+t^2)(s-t)x^2+s^2(s+t)y^2+t^3z^2 }
$$
on the scroll.
The base equations are formed according to Thm.~\labthm.
One computes that indeed we have a complete intersection,
which is non-singular.
\qed

It is very difficult to find solutions to such equations,
and I have not succeeded to do so in 
the specific example. Note that the absence of mixed terms in $x$, $y$ and $z$
on the scroll means that the 
automorphism group of the curve has order at least eight and it operates
on the base space:
given one solution one  finds three
other ones by multiplying all $\xi_i$ or all $\zeta_i$ by $-1$.

\roep Remark \num. 
Alternatively one can start with a $K3$ surface and take a 
general hyperplane section. Therefore we look at complete intersections
of two surfaces of type $2H-5R$ on a scroll of type $(3,3,3,3)$.
Such a $K3$ surface can have infinitesimal deformations
of negative weight (which are always obstructed). 
The lifting matrix for the $K3$ has size $8\times 8$.
The equations $P$ and $Q$ on the scroll are pencils of quadrics.
In general such a pencil has $4$ singular fibres and by taking a suitable
linear combination we may suppose that $P$ has the form
$$
sX^2+tY^2+(s+t)Z^2+(s-t)W^2\;.
$$
The polynomial $Q$ is then a general pencil with
$20$ coefficients, of which one can be made to vanish 
by subtracting a multiple of $P$. This shows that these $K3$ surfaces
depend on $18$ moduli.
Let $Q=(a_{11}s+b_{11}t)X^2+2(a_{12}s+b_{12}t)XY+\dots + (a_{44}s+b_{44}t)W^2$.
Then the lifting matrix is
$$
2\pmatrix{
1   &  0  & 0  &  0 & 0  &  0 & 0  &  0 \cr
0   &  0  & 0  &  1 & 0  &  0 & 0  &  0 \cr
0   &  0  & 0  &  0 & 1  &  1 & 0  &  0 \cr
0   &  0  & 0  &  0 & 0  &  0 & 1  &  -1 \cr
a_{11} &b_{11} &a_{12} &b_{12} & a_{13} &b_{13} & a_{14} &b_{14}  \cr
a_{12} &b_{12} &a_{22} &b_{22} & a_{23} &b_{23} & a_{24} &b_{24}  \cr
a_{13} &b_{13} &a_{23} &b_{23} & a_{33} &b_{33} & a_{34} &b_{34}  \cr
a_{14} &b_{14} &a_{24} &b_{24} & a_{34} &b_{34} & a_{44} &b_{44}  \cr
}\;.
$$
For nonsingular $K3$ surfaces this matrix has at least rank $5$, and it
is possible to write down examples with exactly rank $5$. Rank $4$
can be realised with surfaces with isolated singularities.
An explicit example (with a slightly different basis for the pencil) is
$$
\eqalign
{P&=sX^2+tY^2+(s+t)Z^2\cr
Q&=sX^2-tY^2+(s-t)W^2}
$$
with ordinary double points at $sX=tY=(s+t)Z=(s-t)W=0$.
The hyperplane section $X_3+Z_2+W_1+Y_0=t^3X+s^2tZ+st^2W+s^3Y$
does not pass through the singular points and defines a smooth
tetragonal curve with $e_i=4$. The base space for this curve is
still a complete intersection, but the line corresponding to the 
singular $K3$ surface is a multiple solution.
\endroep

\roep {\rm(\num)} The case $g=16$.
The curves lying on a $K3$ form a codimension one subspace
in the moduli space of tetragonal curves of genus $g=16$.
In terms of the coefficients of the equations of the scroll
one gets an equation of high degree. It makes no sense
to write it. We will not study the most general case $(5,4,4; 6,5)$
but $(5,5,3; 6,5)$. 
These curves form a codimension two subspace in moduli. 
The computations will show that the condition of being 
a hyperplane section has again codimension one.
The lifting matrix need
not have full rank. We have $b_1=2e_3$, and the $g^1_4$
can be composed. 

Suppose that the coefficient of $z^2$ in the first equation on the scroll 
does not vanish. With a coordinate transformation 
we may assume that the equation  has the form $z^2+P_4(s,t;x,y)$
with $P_{4}$ of degree $4$ in $(s\cn t)$ and quadratic in $(x\cn y)$.
Then we can take $Q$ to be without $z^2$ term.
Let $q_{0;1}s^3+\dots+q_{3;1}t^3$ be the coefficient of $xz$
and $q_{0;2}s^3+\dots+q_{3;2}t^3$ that of $yz$.
The rows of the $3\times8$ lifting matrix come only from the
monomial $z$:
$$
\pmatrix{
0  &0  &0  &0   &| &0  &0  &0  &0  &|&1  &0  \cr
0  &0  &0  &0   &| &0  &0  &0  &0  &|&0  &1  \cr
\noalign{\smallskip\hrule}
\noalign{\smallskip}
q_{0;1}  & q_{1;1} &q_{2;1} &  q_{3;1} &| 
&q_{0;2}  & q_{1;2} &q_{2;2} &  q_{3;2} &|&0  &0  \cr
}
$$
The matrix has rank $3$ if some $q_{i,j}$ does not
vanish, but rank $2$ if they all vanish; then the 
surface $\{Q=0\}$ has a singular line.

The deformation variables $\ze_1$, $\ze_2$ vanish.
We have two pure rolling factors deformations $\rho_1$ and $\rho_2$
in the second set of additional equations, and there are scrollar
deformations $\xi_1$, \dots, $\eta_4$. Between those exist a linear
relation given by the third line of the matrix. 
The equations for the base can be written down independently of
this linear relation, because the $\ze_i$ vanish.

We give a specific example:
$z^2+t^4y^2+(s^4+t^4)x^2$ and
$(s^5+t^5)x^2+(s^5-t^5)y^2+q_1(s,t)xz+q_2(s,t)yz$.  
We get the 
following nine equations:
$$
\displaylines{
-2\xi_2\xi_3-2\xi_1\xi_4-2\eta_2\eta_3-2\eta_1\eta_4  \cr 
\xi_1^2- \xi_3^2-2\xi_2\xi_4- \eta_3^2-2\eta_2\eta_4  \cr
2\xi_1\xi_2-2\xi_3\xi_4 -2\eta_3\eta_4 \cr
2\xi_1\xi_3+\xi_2^2- \xi_4^2 - \eta_4^2\cr
2\xi_1\xi_4+2\xi_2\xi_3\cr
\rho_1\xi_1+\rho_2\eta_1-\xi_3^2-2\xi_2\xi_4+\eta_3^2+2\eta_2\eta_4\cr
\rho_1\xi_2+\rho_2\eta_2+\xi_1^2+\eta_1^2-2\xi_3\xi_4+2\eta_3\eta_4\cr
\rho_1\xi_3+\rho_2\eta_3+2\xi_1\xi_2+2\eta_1\eta_2 -\xi_4^2   +\eta_4^2\cr
\rho_1\xi_4+\rho_2\eta_4+2\xi_1\xi_3+\xi_2^2+2\eta_1\eta_3+\eta_2^2\cr
}
$$

Also in general we have $5$ equations $\pi_m$ and $4$ equations
$\rho_1\xi_{m}+\rho_2\eta_m+\chi_m$.
The pure rolling factors equations are never obstructed.
We have as solution to the  equations therefore  the 
$(\rho_1,\rho_2)$-plane with a non reduced structure.
Given a general value of  $(\rho_1,\rho_2)$ we can eliminate
say the $\eta_i$ variables. We are then left with $5$ equations
$\pi_i$ depending only on the $x_i$. Their quadratic parts satisfy
a relation with constant coefficients, but even more is true:
this relation can be  lifted to the equations themselves. So the 
component has multiplicity $16$.
The general fibre over the reduced component has a 
simple elliptic singularity of degree $10$.

To find the other solutions we eliminate $\rho_1$ and $\rho_2$.
This gives the condition
$$
\Rank \pmatrix
{ \chi_1 & \chi_2 & \chi_3 & \chi_4 \cr
  \xi_1 & \xi_2 & \xi_3 & \xi_4 \cr
  \eta_1 & \eta_2 & \eta_3 & \eta_4 \cr}
 \leq 2
$$
which defines a codimension $2$ variety of degree $11$. In general
the $5$ equations $\pi_m$ cut out a subset of codimension
$7$ and degree $352$. But if
$$
\Rank \pmatrix
{
  \xi_1 & \xi_2 & \xi_3 & \xi_4 \cr
  \eta_1 & \eta_2 & \eta_3 & \eta_4 \cr}
 \leq 1
\eqno(R)
$$
the full equations have only solutions in the 
$(\rho_1,\rho_2)$-plane. 
Even if this rank condition defines a codimension $3$ subspace, 
there are always solutions. To see this  we set
$\xi_i=s^{4-i}t^{i-1}\xi$ and $\eta_i=s^{4-i}t^{i-1}\eta$.
The equations  $\pi_m$ are satisfied if
${\partial \over \partial s}P_4(s,t;\xi,\eta)=0$ and
${\partial \over \partial t}P_4(s,t;\xi,\eta)=0$.
This is the intersection of two curves of type $(2,3)$ on 
the scroll $S_{3,3}\cong \P^1\times \P^1$ and there are 
$12$ such intersection points. Those points give multiple solutions.
One can compute that the rank of the Jacobi matrix of the system of 
equations $(R)$
together with the $\pi_m$ is five. By taking a suitable general
example one finds that the multiplicity is in fact $4$, and 
$48$ is the degree of the solution of the system. 

\proclaim Proposition \num.
The general tetragonal cone with invariants $(5,5,3;6,5)$,
which is composed with an involution of genus $4$,
has $302$ smoothing components.
The base space of a non-composed cone can be identified with a 
hyperplane section of the base of the corresponding composed one
and only the smoothing components of lying in this hyperplane
give smoothing components of the non-composed cone.

This means that for  fixed polynomials $P$, $Q$ the existence 
of smoothings depends on an equation of degree $302$
in the eight variables $q_{i;j}$.
For special values the number of smoothing components may go down.
This happens in the specific example given, where the condition
$(R)$ gives two-dimensional `false' solutions.
Here there are only $238$ smoothing components.
Besides the hyperelliptic involution the curve has another
automorphism which acts on the base space. The only solutions I have found 
are easy to see: 
$$
\eta_1=\eta_2=\xi_3=\xi_4=
\xi_1+\xi_2=
\eta_3+\eta_4=
\rho_1+\xi_1=\rho_2+\eta_4=\xi_1^2-\eta_4^2=0
$$
We take $\xi_1=\eta_3=\rho_2=\de$ and $\xi_2=\eta_4=\rho_1=-\de$.
The total space is a surface on a scroll of type $(4,3,3,3)$
with bihomogeneous coordinates $(W,X,Y,Z;s,t)$.
We set $Y_i=y_i$, $X_i=x_{i+2}$ and $Z_i=z_i$  for $i=0,\dots,3$.
The hyperplane section is
$\de = W_2+X_0+X_1+ Y_2 +Y_3$,
so if $\de=0$ we have $X=t^2x$, $Y=s^2y$, $Z=z$ and
$W=-(s+t)(x+y)$.
The lifting equation is now
$q_{0;1}-q_{1;1}+q_{2;2}-q_{3;3}=0$.
One computes that the surface is 
given by
$$
\displaylines{
2X^2+Y^2-2(X-Y)W(s-t)+W^2(s-t)^2+Z^2 \cr
\quad 2Y^2s -XW(s^2-2st+2t^2)+YW(2s^2-2st+t^2) 
    +W^2(s-t)(s^2-st+t^2) \hfill\cr\hfill{}
-XZ(sq_{2;2}+tq_{3;2})-YZ(sq_{0;1}+tq_{1;1}) 
-ZW(s^2q_{0;1}+st(q_{2;2}-q_{3;2})+t^2q_{3;2}) 
\quad
}
$$
This is a $K3$ surface with an $A_1$-singularity.

For even more special values of the coefficients  there
may be higher dimensional smoothing components. This happens e.g. for 
$P=z^2+t^4y^2+s^4x^2$ and the same $Q$ as above,
where the equations $\pi_m$ have the solution
$\xi_1=\xi_2=\eta_3=\eta_4=0$, giving rise to an extra component of degree 15,
which is the cone over three rational normal curves of degree 5.
Then all tetragonal on $Y$ have smoothings, but depending
on the position of the hyperplane the number may increase.

\roep {\rm(\num)} The case $(b_1,b_2)=(8,4)$.
In this case there exist five families of $K3$-surfaces,
three of which have the maximal dimension $18$. The general hyperplane
section of the scroll $S_{8,4,2,0}$ is a scroll  $S_{8,4,2}$
while for both $S_{5,4,3,2}$ and $S_{4,4,4,2}$ it is $S_{5,5,4}$.
One computes that the tetragonal curves of type $(2H-8R,2H-4R)$
on $S_{8,4,2}$ depend on $29$ moduli and those on $S_{5,5,4}$
depend on $34$ moduli.

\proclaim Proposition \num.
The general tetragonal curve of  type $(8,4,2;8,4)$ has only
pure rolling factors extensions. If the $g^1_4$ is composed with an
involution of genus $3$, then there are in general $91$ smoothing
components not of this type.

\roep Remark \num.
The tetragonal curve can be a special hyperplane section of
a $K3$ surface on $S_{7,4,2,1}$, $S_{6,4,2,2}$, $S_{5,4,3,2}$ or 
$S_{4,4,4,2}$. Therefore the genericity assumption cannot be dropped.

\roep Proof.
After a coordinate we may assume that $P$ has the form
$p_8x^2+y^2+p_2xz$. The $g^1_4$ is composed with an
involution of genus $3$ if and only if $p_2\equiv 0$.
In that case $Q$ may be taken in  the form
$q_{12}x^2+q_8xy+z^2$. That the curve is nonsingular implies
that $p_8$ has no multiple roots.
If the $g^1_4$ is not composed, the term $z^2$ may be absent in $Q$,
and $p_8$ may have multiple roots. For the general curve this does not
occur.
We look therefore at curves given by
$$
\eqalign{
P&\colon p_8x^2+y^2+p_2xz\cr
Q&\colon q_{12}x^2+q_8xy+z^2\;.
}
$$
The lifting matrix is a block matrix
$$
\pmatrix {
0 & 2I &0 \cr
\Pi & 0 & 0 \cr
0&0&2I}
$$
with $\Pi$ giving the equations $p_{2,0}\xi_i+p_{2,1}\xi_{i+1}
+p_{2,2}\xi_{i+2}=0$.
There is one pure rolling factors deformation for the first equation,
and $5+1$ for the second. The equation $P$ leads to $7$ base equations
$\pi_m$ in the $8$ variables $\rho$, $\xi_1$, \dots, $\xi_7$.
The $128$ solutions are described above.
The equations coming from $Q$ are
$$
\displaylines {
\rho_1\xi_1+\rho_2\xi_2+\rho_3\xi_3+\rho_4\xi_4+\rho_5\xi_5+\chi_1=0\cr
\rho_1\xi_2+\rho_2\xi_3+\rho_3\xi_4+\rho_4\xi_5+\rho_5\xi_6+\chi_2=0\cr
\rho_1\xi_3+\rho_2\xi_4+\rho_3\xi_5+\rho_4\xi_6+\rho_5\xi_7+\chi_3=0}
$$
We view this as inhomogeneous linear equations for the $\rho_i$.
The coefficient matrix
$$
M=\pmatrix {
\xi_1 & \xi_2 & \xi_3 & \xi_4 & \xi_5\cr
\xi_2 & \xi_3 & \xi_4 & \xi_5 & \xi_6\cr
\xi_3 & \xi_4 & \xi_5 & \xi_6 & \xi_7}
$$
is the transpose of the coefficient matrix of the equations
$p_{2,0}\xi_i+p_{2,1}\xi_{i+1}
+p_{2,2}\xi_{i+2}=0$, viewed as equations for the coefficients of $p_2$.
If for a given solution of the equations $\pi_m$ the matrix 
$M$ has not full rank, then  there exists 
a non-composed pencil admitting the same solution.
But then also $p_{2,0}\chi_1+p_{2,1}\chi_{2}
+p_{2,2}\chi_{3}=0$, an equation which in general is not satisfied.
We have $8$ solutions which lie on a rational normal curve and $28$
solutions on the secant variety of this curve. The equations of the
secant variety are the maximal minors of $M$. Only for $91$ solutions
the matrix $M$ has full rank.
\qed

In the general case we get  components of dimension $3+1$
(the $y$-rolling factors deformation does not enter the equations), 
for solutions not on the rational curve, but on its secant variety
the component has dimension $5$, while we get a $6$-dimensional component if
$p_8$ and $q_{12}$ have a common root. This does not contradict
the fact that all smoothing components of Gorenstein surface singularities
have the same dimension, because we here only look at the restriction to
negative degree.

\proclaim Proposition \num.
The general hyperplane section of a $K3$  of type
$(5,4,3,2;8,4)$ or $(4,4,4,2;8,4)$ is a tetragonal curve
of type $(5,5,4;8,4)$, which lies on a rational surface with
two double points.

\roep Proof.
We use coordinate transformations on the scroll to bring the hyperplane
section into a normal form, while we suppose 
the coefficients of the equations to be general.
Let as usual $(X,Y,Z,W;s,t)$ be coordinates on the scroll.
Let the hyperplane section be $a_0W_0+a_1W_1+a_2W_2+ \cdots=
(a_0s^2+a_1st+a_2t^2)W+\cdots$. By a transformation in $(s,t)$
we achieve that $a_0=a_2=0$, so the equation is $W_1+\cdots$.
First consider the case $(4,4,4,2)$. By a suitable transformation
$w\mapsto W+a_2(s,t)X+b_2(s,t)Y+c_2(s,t)Z$ we remove all terms
with index $1$, $2$ or $3$, leaving $W_1+a_0X_0+b_0Y_0+c_0Z_0
+a_4X_4+b_4Y_4+c_4Z_4$. Taking $a_0X+b_0Y+c_0Z$ as new $X$ and
$a_4X+b_4Y+c_4Z$ as new $Y$ brings us finally to $X_4+W_1+Y_0$.
With coordinates $(x,y,z;s,t)$ for the scroll $S_{5,5,4}$ we
get the hyperplane section by setting $Z=z$, $X=sx$, $Y=ty$
and $W=-t^3x-s^3y$. The equation $P$ does not involve the variable
$W$ so we have quadratic singularities if $sx=ty=z=0$, which gives
the points $s=y=z=0$ and $t=x=z=0$.

In the case $(5,4,3,2)$ we can achieve $W_1+Z_0+Z_3$ and we
get the curve by $X=x$, $Y=z$, $Z=sty$ and $W=-(s^3+t^3)y$.
The equation $P\colon p_2X^2+p_1XY+p_0Y^2+XZ$ now gives
$p_2x^2+p_1xz+p_0z^2+stxy$, which for general $p_i$ has 
singular points at $x=z=st=0$.
\qed

To investigate the sufficiency of these conditions we look
at the general cone of type $(5,5,4;8,4)$.
We may suppose that  $P$ has the form $z^2+P_2(x,y)$. 
The equation $P_2(x,y)$ describes a curve of type $(2,2)$ on $S_{5,5}\cong 
\P^1 \times \P^1$. If this curve has a singular point, we may
assume that it lies in the point $x=s=0$. Under the assumption 
that the coefficient of $stxy$ does not vanish we can transform the
equation  into the form $(as^2+bt^2)x^2+2stxy+cs^2y^2$ and
unfolding the singularity we get the equation
$$
P_2= (as^2+bt^2)x^2+2stxy+(cs^2+dt^2)y^2\;.
$$
One can then write out the lifting conditions and base equations
coming from the equation $P$. The result is that they have only
trivial solutions if and only if
$abcd((ad+bc-1)^2-4abcd)\neq 0$, if and only if the curve $P_2$
is nonsingular.
If a singularity is present we assume it to be in $x=s=0$,
so $d=0$. The equation $Q$ gives three base equations, in which
$2+2$ pure rolling factors variables can enter.
We analyse what happens if there is a second singularity.
For $b=d=0$ the equation $P_2$ is divisible by $s$, and we do not find
extensions.
In case $a=d=0$ the curve $P_2$ splits into two curves of type $(1,1)$;
we get two components with deformed scroll $S_{4,4,4,2}$.
For $c=d=0$ we have intersection
of a line with a curve of type $(2,1)$ and we find  two components
with deformed scroll $S_{5,4,3,2}$.

\comment 
To investigate the sufficiency of these conditions we write 
base equations for a tetragonal cone of type $(5,5,4;8,4)$.
The equation $P$ need not contain the term $z^2$. We look at the general
case and suppose that  $P$ has the form $z^2+P_2(x,y)$. The equation
$P_2(x,y)$ describes a curve of type $(2,2)$ on $S_{5,5}\cong 
\P^1 \times \P^1$. If this curve has a singular point, we may
assume that it lies in the point $x=s=0$. Under the assumption 
that the coefficient of $stxy$ does not vanish we can transform the
equation  into the form $(as^2+bt^2)x^2+2stxy+cs^2y^2$ and
unfolding the singularity we get the equation
$$
P_2= (as^2+bt^2)x^2+2stxy+(cs^2+dt^2)y^2\;.
$$
This curve has a singularity if and only if
$abcd((ad+bc-1)^2-4abcd)=0$. Because our equation is invariant under the
reflection $(x\cn y; s\cn t)\mapsto (-x\cn y; -s\cn t)$, the curve
is has two singular points if $(ad+bc-1)^2=4abcd$. It splits then
into two curves of type $(1,1)$. The same happens if $a=d=0$. 
If $b=d=0$, the equation is divisible by $s$, so the curve
splits into a line and a curve of type $(1,2)$. The intersection
of a line with a curve of type $(2,1)$ occurs for $c=d=0$.

Because of the presence of the term $z^2$ in $P$ the lifting 
conditions force all $\ze_i$ to vanish.
We write the lifting matrix (divided by $2$) for $P_2$:
$$
\pmatrix
{ a & 0 & b & 0 & 0 & 1 & 0 & 0 \cr
  0 & a & 0 & b & 0 & 0 & 1 & 0  \cr
  0 & 1 & 0 & 0 & c & 0 & d & 0 \cr
  0 & 0 & 1 & 0 & 0 & c & 0 & d 
}
$$
Using the lifting conditions we can write the seven base equations
as follows.
$$
\displaylines{
b\xi_1\xi_2+\xi_1\eta_1+d\eta_1\eta_2 \cr
a\xi_1^2-b\xi_2^2+c\eta_1^2-d\eta_2^2 \cr
a\xi_1\xi_2+\xi_2\eta_2+c\eta_1\eta_2 \cr
a\xi_2^2-b\xi_3^2+c\eta_2^2-d\eta_3^2 \cr
a\xi_2\xi_3+\xi_3\eta_3+c\eta_2\eta_3 \cr
a\xi_3^2-b\xi_4^2+c\eta_3^2-d\eta_4^2 \cr
a\xi_3\xi_4+\xi_4\eta_4+c\eta_3\eta_4 }
$$
If the discriminant $abcd((ad+bc-1)^2-4abcd)$ does not vanish, these
equations (together with the lifting equations) 
have only the trivial solution.
Otherwise we may assume that there is a singular point in $x=s=0$,
so $d=0$. Then we find solutions if there are at least two
singular points. In case $b=0$, but $ac\neq0$, we find $\xi_1=\xi_2=\xi_3=
\eta_1=\eta_2=\eta_3=0$ and $\xi_4\eta_4=0$.
The equation $Q$ gives three base equations, in which
$2+2$ pure rolling factors variables can enter. If only
$\xi_4$ is non-zero, the third equation contains the term
$\xi_4\rho_2$. The first equation vanishes, while the second reduces to
$-q_6\xi_4^2$, where $q_6$ is the coefficient of $t^2x^2$ in $Q$.
Nonsingularity of the tetragonal curve gives 
the condition that $Q$ does not pass through the singular point
$z=y=s=0$, so $q_6\neq0$. In the same way one finds a
contradiction if $\eta_4\neq0$.  

If $c=d=0$ we get $\xi_1=\xi_2=\xi_3=\xi_4=\eta_2=\eta_3=0$.
In this case the equations reduce to 
$$
\matrix{
\eta_1 \sigma_1&-&q'_1\eta_1^2& \cr
&&q'_0\eta_1^2&-&q'_6\eta_4^2 \cr
\eta_4 \sigma_2&&&+&q'_5\eta_4^2
}
$$
with the $q'_i$ the coefficients in the $y^2$-term.
Nonsingularity of the tetragonal curve gives that $q_0q_6\neq0$,
so we get two components. In both cases $\eta_4$ is a multiple
of $\eta_1$, so the deformed scroll is $S_{5,4,3,2}$.

Other solutions of the base equations occur for 
intersection of two curves of type $(1,1)$. If $d=bc-1$
they intersect in one point with the same tangent.
The general case occurs for $a=d=0$.
Then $\xi_2=\xi_3=\xi_4=\eta_1=\eta_2=\eta_3=0$.
Now the second equation reduces to $q_0\xi_1^2-q'_6\eta_4^2$
and we get two components with deformed scroll $S_{4,4,4,2}$.
\endcomment

\roep Remark \num.
For the general tetragonal cone with large $g$ we found
$\dim T^1(-1)=9$, but all deformations are obstructed.
For special curves extensions may exist; 
also the dimension can be higher. Both conditions seem to be 
independent. As the number of base equations we find is always $g-7$,
having more variables increases the chances of finding solutions.
In the borderline case studied above this may suffice to force the
existence, but in general it does not.  On the other, taking
a general hyperplane section of a general tetragonal $K3$ surface will
give a cone with $\dim T^1(-1)=9$. It would be interesting to
find a property of a canonical curve which gives a sufficient
condition for the existence of an extension.

\vskip 2cm

\roep {\bf Acknowledgement}. I thank Jim Brawner for sending me his preprint.
The computer computations were done with the program
{\sl Macaulay\/} \cite\rfbs. 

\vfil

\beginsection References.

\frenchspacing
\parindent=0pt
\def\item{\par\hangindent=20pt}

\item
        Dave Bayer and Mike Stillman, 
       {\sl Macaulay: A system for computation in
        algebraic geometry and commutative algebra.\/}
        Computer software available via anonymous ftp from
       {\tt //math.harvard.edu/Macaulay/}.

\item 
   James N. Brawner,
   {\sl Tetragonal Curves, Scrolls and $K3$ Surfaces\/}.
   Trans. Amer. Math. Soc. {\bf 349} (1997), 3075--3091.
   
\item 
   James N. Brawner,
   {\sl The Gaussian-Wahl Map for Tetragonal Curves\/}.
   Preprint 1996.

\item 
  Ciro Ciliberto and Rick Miranda,
  {\sl Gaussian maps for certain families of canonical curves\/}.
  In: Complex Projective Geometry. London Math. Soc. Lecture Note Ser.
  {\bf  179}, Cambridge Univ. Press, Cambridge, 1992, pp. 106--127.

\item
   Ciro  Ciliberto and Guiseppe Pareschi,
   {\sl Pencils of minimal degree on curves on a $K3$ surface\/}. 
   J. Reine Angew. Math. {\bf 460} (1995), 15--36. 

\item
   Ron Donagi and David R. Morrison,
   {\sl Linear systems on $K3$-sections\/}. 
   J. Differential Geom. {\bf 29} (1989), 49--64. 

\item 
   R. Drewes and J. Stevens,
   {\sl Deformations of Cones over
   Canonical Trigonal Curves\/}.
   Abh. Math. Sem. Univ. Hamburg {\bf 66} (1996), 289--315.

\item  
  Patrick Du Val, 
  {\sl On rational surfaces whose prime sections are canonical curves\/}. 
  Proc. London Math. Soc. (2) {\bf 35}, (1933), 1--13.

\item
  D. H. J. Epema, 
  {\sl Surfaces with Canonical Hyperplane Sections\/}.
  Thesis, Leiden 1983. Also:  CWI Tract {\bf 1}, 
   Centrum voor Wiskunde en Informatica, Amsterdam, 1984.

\item 
  Miles Reid, 
  {\sl Surfaces with $p_g=3, K^2=4$ according to E. Horikawa 
  and D. Dicks\/}. 
  Text of a lecture, Univ. of Utah  and Univ. of Tokyo 1989.

\item 
  Miles Reid, 
  {\sl Chapters on Algebraic Surfaces\/}. 
  In: Complex algebraic geometry (Park City 1993), 
  IAS/Park City Math. Ser. {\bf3}, 
  Amer. Math. Soc.,  1997, pp. 3--159.

\item 
  Karl Rohn,
  {\sl Ueber die Fl\"achen vierter Ordnung mit dreifachem Punkte\/}.
  Math. Ann. {\bf 24} (1884), 55--151.

\item 
  Frank-Olaf Schreyer,
  {\sl Syzygies of Canonical Curves and Special Linear
  Series\/}. 
  Math. Ann. {\bf 275} (1986), 105--137.

\item 
  Jan Stevens, 
  {\sl Deformations of cones over hyperelliptic curves\/}.  
  J. Reine Angew. Math. {\bf 473} (1996), 87--120.

\item 
  Jonathan  Wahl,
  {\sl On the cohomology of the square of an ideal sheaf\/}. 
  J. Algebraic Geom. {\bf 6} (1997) 481--511.

\item 
  Jonathan Wahl,
  {\sl Hyperplane sections of Calabi-Yau varieties\/}. 
  Preprint 1998.

\vfill

\parindent=0pt
\parskip=0pt
{\obeylines
Address of the author:
Jan Stevens
Matematik
G\"oteborgs universitet
Chalmers tekniska h\"ogskola
SE 412 96 G\"oteborg, Sweden
e-mail: stevens@math.chalmers.se
}

\bye